\documentclass[titlepage,11pt]{article}

\usepackage{amssymb,latexsym,amsmath,amsthm,verbatim}
\numberwithin{equation}{section}

\newcommand{\eps}{\varepsilon}
\newcommand{\be}{\begin{equation}}
\newcommand{\ee}{\end{equation}}
\newcommand{\ba}{\begin{array}}
\newcommand{\ea}{\end{array}}
\newcommand{\ds}{\displaystyle}

\newcommand{\R}{\mathbb{R}}
\newcommand{\C}{\mathbb{C}}
\newcommand{\X}{\mathcal{X}}
\newcommand{\Y}{\mathcal{Y}}
\newcommand{\ZZ}{\mathcal{Z}}
\newcommand{\Z}{\mathbb{Z}}
\newcommand{\N}{\mathbb{N}}
\newcommand{\sgn}{\mathop{{\rm sgn}}}
\newcommand{\re}{\mathop{{\rm Re}}}
\newcommand{\col}{\mathop{{\rm col}}}

\newcommand{\codim}{\mathop{{\rm codim}}}
\newcommand{\ran}{\mathop{{\rm ran}}}
\newcommand{\LL}{\mathcal{L}}

\newcommand{\diag}{\mathop{{\rm diag}}}
\newcommand{\bxx}{\begin{it}}
\newcommand{\exx}{\end{it}}
\newcommand{\bb}{\hspace*{-.08in}}
\newcommand{\vv}{\vspace{.2in}}
\newcommand{\strt}{\rule{0mm}{5mm}}

\begin{document}
\baselineskip 1.6em
\title{Universality of Crystallographic Pinning}

\author{John Mallet-Paret\thanks{Partially supported by NSF DMS-0500674}
\\Division of Applied Mathematics
\\Brown University
\\Providence, RI 02912
\and
Aaron Hoffman\thanks{Partially supported by NSF DMS-0603589}
\\Department of Mathematics and Statistics
\\Boston University
\\Boston, MA 02215
\\$\;$}

\date{October 31, 2008}

\maketitle

\begin{abstract}
We study traveling waves for reaction diffusion equations on the spatially discrete domain $\Z^2$. 
The phenomenon of crystallographic pinning occurs when traveling waves become pinned in certain directions
despite moving with non-zero wave speed in nearby directions.  In \cite{JMP-CP} it was shown that crystallographic
pinning occurs for all rational directions, so long as the nonlinearity is close to the sawtooth. 
In this paper we show that crystallographic pinning holds in the horizontal and vertical directions for bistable
nonlinearities which satisfy a specific computable generic condition. 
The proof is based on dynamical systems.  In particular,
it relies on an examination of the heteroclinic chains which occur as singular limits of wave
profiles on the boundary of the pinning region.
\end{abstract}

\section{Introduction}
The setting for this paper is traveling waves for lattice differential equations of reaction diffusion type.
A lattice differential equation (LDE) is an infinite system of coupled ordinary differential equations,
where each ODE represents the dynamics at a single point on a spatial lattice.  A simple LDE is
\begin{equation}
\dot{u}_i = d(u_{i+1} + u_{i-1} - 2u_i) - f(u_i) \label{1DLDE}
\end{equation}
where $i \in \Z$ is a spatial index, $d \in \R$ is the coupling constant,
$f: \R \to \R$ is a given function, and each $u_i$ is a function of a single variable $t$.
If $d > 0$ then equation \eqref{1DLDE} may be regarded as a reaction-diffusion equation:
The first term $d(u_{i+1} + u_{i-1} - 2u_i)$ is a discrete second derivative which provides the diffusion,
while the second term $f(u_i)$ is the reaction term.

If we denote $d = \frac{1}{h^2}$, then equation \eqref{1DLDE} is obtained from the PDE
\begin{equation}
u_t = u_{xx} - f(u) \label{1DPDE}
\end{equation}
upon replacing the term $u_{xx}$ by a standard central difference approximation with grid size $h$.
The limit $h \to 0$, that is $d \to \infty$ in equation \eqref{1DLDE}, corresponds at least formally to equation \eqref{1DPDE}.
In this study we are interested in spatially discrete systems such as
\eqref{1DLDE} which are not necessarily close to the PDE limit in this sense,
namely $d$ need not be large, although we do assume $d>0$.
In fact, without loss we may take $d=1$ by rescaling time and redefining $f$.
This normalization serves to emphasize that the grid is not particularly small and is not meant to approximate a continuum.
Thus we study the equation
\begin{equation}
\dot{u}_i = u_{i+1} + u_{i-1} - 2u_i - f(u_i) \label{1DLDEnorm}
\end{equation}
and its higher-dimensional analogs in what follows.

Lattice models are widely used in applications such as solid state physics, materials science, and physiology;
see \cite{Bell,Erneux,Firth,JMP-LNM,Perez} and the references therein.

We assume the nonlinearity $f$ is of {\bf bistable} type, in particular that 
\be
\ba{l}
f(\pm 1)=0,\qquad f(a)=0,\qquad f'(\pm 1)>0,\qquad f'(a)<0,\\
\\
f(u)>0\hbox{ for }u\in(-1,a)\cup(1,\infty),\\
\\
f(u)<0\hbox{ for }u\in(-\infty,-1)\cup(a,1),
\ea
\label{bistable}
\ee
for some $a \in (-1,1)$.
Generally, only the values of $u$ for $|u|\le 1$ will be relevant
for our arguments, but in several places (in particular some arguments by contradiction)
it will be convenient to assume that $f(u)$ is smoothly extended for $|u|>1$, with
the sign condition as in \eqref{bistable}.
We in fact take a family of bistable functions $f=f(u,a)$ parameterized by $a\in(-1,1)$,
in addition satisfying the monotonicity condition
\be
\frac{\partial f}{\partial a}(u,a) > 0
\hbox{ for }u \in(-1,1)\hbox{ and }a \in (-1,1)
\label{normal}
\ee
in $a$. To be precise, let us define a set $\mathcal{N}$ of functions by
$$
\ba{lcl}
\mathcal{N} &\bb = &\bb \{f:[-1,1]\times(-1,1)\to\R\;|\;f(\cdot,\cdot)\hbox{ is }C^{2}\hbox{ smooth, with }f(\cdot,a)\\
\\
&\bb &\bb \hbox{satisfying }\eqref{bistable}\hbox{ for every }(u,a)\in[-1,1]\times(-1,1),\hbox{ and }\eqref{normal}\hbox{ holding}\}.
\ea
$$
We say that $f:[-1,1]\times(-1,1)\to\R$ is a {\bf normal family} if $f\in\mathcal{N}$.
We do not endow the set $\mathcal{N}$ with a topology, although we shall do so with
certain subsets of $\mathcal{N}$.
The parameter $a\in(-1,1)$ is known as the {\bf detuning parameter}.
The function
$$f(u,a)=(u^2-1)(u-a)$$
furnishes a simple example of a normal family.

The bistable reaction term $f$ in equation \eqref{1DLDEnorm} pushes the system toward spatial heterogeneity;
it forces $u_i$ toward $+1$ when $u_i > a$ and toward $-1$ when $u_i < a$.
By constrast, the diffusion term $u_{i+1} + u_{i-1} - 2u_i$ promotes
spatial homogeneity, forcing $u_i$ toward the average $\frac{1}{2}(u_{i+1} + u_{i-1})$ of its neighbors.
This competition between reaction and diffusion,
spatial heterogeneity and spatial homogeneity,
is what gives reaction diffusion equations of bistable type the richness to
support both spatially chaotic patterns and traveling waves \cite{JMP-CP}.

By a traveling wave solution of equation \eqref{1DLDEnorm} we mean a solution of the form 
\begin{equation}
u_i(t) = \phi(i-ct),\qquad i \in \Z,
\label{wave}
\end{equation}
for some function $\phi : \R \to \R$, the so-called {\bf wave profile},
and some $c \in \R$, the {\bf wave speed}.  In this paper we consider monotone traveling waves which connect the spatially homogeneous equilibria at $\pm 1$,
that is, $\phi$ satisfies the boundary conditions
\be
\phi(-\infty) = -1,\qquad
\phi(\infty) = 1,
\label{BC}
\ee
and $\phi(\xi)$ is monotone in $\xi$.  

Substitution of \eqref{wave} into \eqref{1DLDEnorm} shows that the
function $\phi$ must satisfy the {\bf wave profile equation}
\begin{equation}
-c\phi'(\xi) = \phi(\xi+1)+\phi(\xi-1)-2\phi(\xi) - f(\phi(\xi),a) \label{WP}
\end{equation}
for $\xi\in\R$ if $c\ne 0$, or for $\xi\in\Z$ if $c=0$.  Conversely, any solution of \eqref{WP} generates a traveling wave solution of \eqref{1DLDEnorm}.  Note that equation \eqref{WP} is a
differential-difference equation if $c\ne 0$, and a difference equation if $c=0$.  This difference in character of the wave profile equation between $c=0$ and $c \ne 0$ is one of the chief reasons that the LDE \eqref{1DLDEnorm} exhibits behavior which is not present in the PDE \eqref{1DPDE}.

In the PDE case \eqref{1DPDE} with $x \in \R$, traveling wave solutions have the form 
$u(t,x) = \phi(x-ct)$ where
$\phi$ satisfies the ordinary differential equation
\begin{equation}
-c\phi' = \phi'' - f(\phi,a).
\label{WPODE}
\end{equation}

For both the LDE \eqref{1DLDEnorm} and PDE \eqref{1DPDE} problems above, it is known that there is a unique wavespeed $c$ at which a monotone traveling wave exists, that is, a monotone solution to either equation \eqref{WP} or \eqref{WPODE}, which satisfies the boundary conditions \eqref{BC}.
Of course the wave speed $c = c(a)$ depends on the parameter $a \in (-1,1)$ as does the wave profile $\phi(\xi) = \phi(\xi;a)$.
Moreover, the wave speed is continuous and nondecreasing in $a$.
In the PDE case $c(a)$ is strictly increasing, in fact $c'(a) > 0$ for all $a\in(-1,1)$.
The LDE case differs from the PDE case in that the wave speed may be zero for an open set of $a$.
In particular, there are quantities
$$-1 \le a_- \le a_+ \le 1$$
such that
$$\ba{l} c(a) = 0\hbox{ for }a \in [a_-,a_+]\cap(-1,1), \\
\\
c(a) > 0\hbox{ for }a \in (a_+,1), \\
\\
c(a) < 0\hbox{ for }a \in (-1,a_-), \ea $$
with $c(a)$ depending smoothly on $a$ and $c'(a) > 0$, both whenever $c(a) \ne 0$.
In case a strict inequality $a_- < a_+$ holds,
we say that {\bf pinning} or {\bf propagation failure} occurs; the wave is pinned and cannot propagate when $a$ is between these values.
The interval $(a_-,a_+)$ is called the {\bf pinning interval} and its length measures
the severity of the pinning.
Pinning was observed by Bell \cite{Bell}, Bell and Cosner \cite{BellCosner}, and Keener \cite{Keener}
and has been studied extensively (see for example \cite{Chmaj,Carpio,Fath,Matthies}).  

To see why pinning can occur in the LDE case, observe that when $c = 0$
equation \eqref{WP} reduces to a pure difference equation 
$$
0 = \phi(\xi+1) + \phi(\xi-1) - 2\phi(\xi) - f(\phi(\xi),a),\qquad \xi \in \Z,
$$
which is equivalent to the discrete-time dynamical system
$$
\ba{lcl}
q_{n+1} &\bb = &\bb 2q_n-r_n+f(q_n,a),\\
\\
r_{n+1} &\bb = &\bb q_n
\ea
$$
under the transformation $r_n = \phi(n)$ and $q_n = \phi(n+1)$. 
A solution of \eqref{WP} satisfying the boundary conditions \eqref{BC} thus
corresponds to a heteroclinic connection between the equilibria $(-1,-1) \in \R^2$ and $(1,1) \in \R^2$ in this dynamical system.
As both these equilibria are saddles, such a heteroclinic connection may lie on a transverse intersection of their stable and unstable manifolds, in which case this connection will persist as $a$ varies in some interval, with $c(a) = 0$ throughout this interval.
The maximal interval on which the connection persists is thus $[a_-,a_+]$.

By contrast, pinning does not occur for the PDE \eqref{1DPDE} because here the traveling wave
corresponds to a saddle-saddle connection in the continuous time planar dynamical system \eqref{WPODE}.
Such heteroclinic solutions do not generally persist under perturbations, and in fact
$c(a) = 0$ for a unique value of $a$.

Our interest is in propagation failure for LDE's as it occurs in higher-dimensional lattices
such as $\Z^k\subseteq\R^k$. The phenomena here are more subtle
because all of the quantities mentioned thus far ---
the wave speed $c$, the wave profile $\phi$, and the pinning interval $(a_-,a_+)$ ---
depend now on the direction of propagation.
We take the lattice $\Z^2$, for which the analog of equation \eqref{1DLDEnorm} is
\begin{equation}
\dot{u}_{i,j} = (\Delta u)_{i,j} - f(u_{i,j},a),
\label{LDE}
\end{equation}
with the nonlinearity $f$ as before.
Here the coordinate $(i,j)$ indexes a spatial point in $\Z^2$, each $u_{i,j}$ is a function of $t$ as before,
and the discrete laplacian $\Delta$ is given by 
\begin{equation} \label{2dlap}
(\Delta u)_{i,j} = u_{i+1,j} + u_{i-1,j} + u_{i,j+1} + u_{i,j-1} - 4u_{i,j}.
\end{equation}
A traveling wave solution of \eqref{LDE} is a solution of the form
$$
u_{i,j}(t) = \phi(i\kappa+j\sigma-ct)
$$
for some $(\kappa,\sigma) \in \R^2 \setminus \{(0,0)\}$,
termed the {\bf direction vector}, and some function $\phi:\R\to\R$.
The wave profile equation satisfied by $\phi$ now takes the form
\begin{equation}
-c\phi'(\xi) = \phi(\xi+\sigma) + \phi(\xi-\sigma) + \phi(\xi+\kappa)+\phi(\xi-\kappa) - 4 \phi(\xi) - f(\phi(\xi),a).
\label{2DWP}
\end{equation}
It is known that for each direction vector $(\kappa,\sigma) \in \R^2 \setminus \{0\}$
and each $a \in (-1,1)$, there is a unique wave speed $c = c(a,(\kappa,\sigma))$
such that equation \eqref{2DWP} admits a monotone solution satisfying the boundary
conditions \eqref{BC}.
Moreover, when $c \ne 0$ this solution $\phi = \phi(\xi;a,(\kappa,\sigma))$ is unique up to translation.
The wave speed $c(a,(\kappa,\sigma))$ depends continuously on $a$ and on $(\kappa,\sigma)$.
For each $(\kappa,\sigma)$ it is nondecreasing in $a$, and as before is smooth in $a$ and satisfies
$\frac{\partial c(a,(\kappa,\sigma))}{\partial a}>0$ when $c(a,(\kappa,\sigma))\ne 0$.
Also as before, we have quantities $a_\pm=a_\pm(\kappa,\sigma)$ characterized by
$$[a_-(\kappa,\sigma),\; a_+(\kappa,\sigma)]\cap(-1,1)
=\{a\in(-1,1)\;|\; c(a,(\kappa,\sigma)) = 0\}.$$
Writing $(\kappa,\sigma)=(r\cos\theta,r\sin\theta)$ for some $r>0$ and $\theta\in\R$,
one easily checks by rescaling the independent variable $\xi$ in \eqref{2DWP}
by a factor $r$ that
$$c(a,(r\cos\theta,r\sin\theta)) = r c(a,(\cos \theta, \sin\theta)),
\qquad a_\pm(r\cos\theta,r\sin\theta) = a_\pm(\cos\theta,\sin\theta).$$
We will sometimes abuse notation by writing $a_\pm(\theta)$ instead of $a_\pm(r\cos \theta, r\sin \theta)$.

The functions $a_\pm(\theta)$ are the central objects of study in this paper.
For definiteness we study $a_+(\theta)$.
Although $c(a,(\kappa,\sigma))$ depends continuously on both $a$ and $(\kappa,\sigma)$,
the function $a_+(\theta)$ need not depend continuously on $\theta$.
However, $a_+(\theta)$ is upper semi-continuous in $\theta$, that is,
\be
\limsup_{\theta \to \theta_0} a_+(\theta)\le a_+(\theta_0)
\label{CP}
\ee
holds for every $\theta_0$. This is an immediate consequence of the continuity
of the function $c(\cdot,\cdot)$.
\vv

\noindent {\bf Definition.}
We say that
{\bf crystallographic pinning} occurs for the system \eqref{LDE} in the direction $\theta_0$
in case either the inequality \eqref{CP} is strict, or the analogous inequality for
$a_-(\theta)$ is strict.
\vv

It has been established in \cite{JMP-CP} that for $f$ sufficiently close to the sawtooth function
$f_0(u,a) = u-\sgn(u-a)$,
crystallographic pinning occurs in every direction $\theta_0$ for which $\tan \theta_0$ is rational.  Numerical studies \cite{Hupkes} suggest that  crystallographic pinning occurs in the vertical and horizontal directions when $f$ is the cubic nonlinearity.

The goal of this paper is to give a specific generic condition on a general nonlinearity $f$
under which crystallographic pinning occurs at $\theta_0 = 0$,
that is, with $(\kappa,\sigma)=(1,0)$. This is given in Theorem 1.1 below.
Before stating that result, let us introduce two conditions which will be needed.
These conditions do not necessarily hold for an arbitrary $f$, but rather will
be taken as hypotheses in our main result.
\vv

\noindent {\bf Condition A\@.}
There exists $p\in\ell^\infty(\Z)$, denoted $p=\{p_n\}_{n\in\Z}$, satisfying
\begin{equation}
p_{n+1} + p_{n-1} - 2p_n = f(p_n,a_+(0)),\qquad n \in \Z,
\label{pdiff}
\end{equation}
and which also satisfies the boundary and monotonicity conditions
\be
\lim_{n \to \pm \infty} p_n = \pm 1,\qquad
p_n \le p_{n+1}\hbox{ for }n\in\Z.
\label{BM}
\ee
Moreover, such $p$ is unique up to a shift in the index $n$.
\vv

\noindent {\bf Condition B\@.}
Condition A holds. Further, if $v\in\ell^\infty(\Z)\setminus\{0\}$ satisfies
\begin{equation}
v_{n+1} + v_{n-1} - 2v_n = f'(p_n,a_+(0))v_n,\qquad n\in\Z, \label{lindiff}
\end{equation}
where $p$ is as in Condition A and where
we denote $f'(u,a)=\frac{\partial f(u,a)}{\partial u}$, then the
quantity
\begin{equation}
B=\frac{1}{2}\sum_{n=-\infty}^\infty f''(p_n,a_+(0))v_n^3
\label{bsum}
\end{equation}
satisfies $B\ne 0$.
\vv

We now state the main results of this paper.
\vv

\noindent {\bf Theorem 1.1.} \bxx
Assume that Condition B holds.
Then the inequality \eqref{CP} is strict at $\theta_0=0$
and so crystallographic pinning occurs in the direction $\theta_0 = 0$.
\exx
\vv

\noindent {\bf Theorem 1.2.} \bxx
Condition $B$ is generic in the following sense.  Fix any $f_0 \in \mathcal{N}$ and define the set
\be
C^{2}_+=\{\gamma\in C^{2}[-1,1]\;|\;\gamma(u)>0\hbox{ for every }u\in[-1,1]\},
\label{C3}
\ee
noting that for every $\gamma\in C^{2}_+$ we have $f\in\mathcal{N}$, where $f(u,a)=\gamma(u)f_0(u,a)$.
Let the set $C^{2}_+$ be endowed with the usual $C^{2}$ topology, and so making it an open subset
of the Banach space $C^{2}[-1,1]$.
Then the set $\mathcal{G}(f_0)\subseteq C^{2}_+$ defined as
$$
\mathcal{G}(f_0)=\{\gamma\in C^{2}_+\;|\;\gamma f_0\hbox{ satisfies Condition B}\}
$$
is a residual subset of $C^{2}_+$.
\exx
\vv

Actually, we need the following two propositions in order to ensure the above conditions are
well-defined (for example, to ensure the convergence of the sum in Condition B).
\vv

\noindent {\bf Proposition 1.3.} \bxx
Assume that $c(a,(1,0)) = 0$, equivalently, that $a\in[a_-(0),a_+(0)]$.
Then there exists $p\in \ell^\infty(\Z)$ satisfying 
\begin{equation}
p_{n+1} + p_{n-1} - 2p_n = f(p_n,a),\qquad n \in \Z,
\label{scalardiff}
\end{equation}
along with \eqref{BM}.
Moreover, any such monotone $p$ is strictly monotone, that is $p_n<p_{n+1}$ for all $n \in \Z$.
\exx
\vv

\noindent {\bf Proof.}
The difference equation \eqref{scalardiff} is nothing more than the wave profile
equation \eqref{2DWP} with $c = 0$ and $\theta = 0$, that is, $(\kappa,\sigma)=(1,0)$.
As $c(a,(1,0)) = 0$, this equation has a  monotone heteroclinic solution joining $\pm 1$.
It only remains to show that this solution is strictly monotone.

To show strict monotonicity, suppose to the contrary that $p_{m+1} = p_m$ for some $m\in\Z$.
Then $f(p_m,a) = p_{m-1} - p_m \le 0$.  But $f(p_m,a) = f(p_{m+1},a) = p_{m+2}-p_{m+1} \ge 0$.
Therefore $f(p_m,a)=f(p_{m+1},a)= 0$ and $p_{m-1} = p_m = p_{m+1} = p_{m+2}$.
Continuing in this fashion, we see that $p_n$ is constant in $n$, which contradicts
the fact that it is a heteroclinic connection between $-1$ and $+1$.
\qed
\vv

Equation \eqref{lindiff} can be expressed as $Lv=0$, with the operator $L\in\mathcal{L}(\ell^\infty(\Z))$
given by
\be
L=S+S^{-1}-2I-f'(p,a_+(0)).
\label{lop}
\ee
Here $S\in\mathcal{L}(\ell^\infty(\Z))$ is the shift operator defined as
$$
(SX)_n=x_{n+1}\hbox{ for }X=\{x_n\}_{n\in\Z}\in\ell^\infty(\Z),
$$
and by a slight abuse of notation $f'(p,a_+(0))\in\mathcal{L}(\ell^\infty(\Z))$
denotes the diagonal operator with entries $f'(p_n,a_+(0))$.
\vv

\noindent {\bf Proposition 1.4.} \bxx
Assume that Condition A holds, with $p$ as stated there, and let
$L\in\mathcal{L}(\ell^\infty(\Z))$ be as in \eqref{lop}.
Then there exists $v\in\ell^\infty(\Z)\setminus\{0\}$ satisfying \eqref{lindiff},
that is, $Lv=0$.
The vector $v$
is unique up to scalar multiple, and thus
$$
\ker(L)=\{av\;|\;a\in\R\}.
$$
Further, $v$ can be chosen to satisfy
$$v_n>0,\qquad n\in\Z,$$
and its coordinates enjoys the estimate
\be 
v_n\le K\mu^{|n|},\qquad n\in\Z,
\label{expest}
\ee
for some $K>0$ and $0<\mu<1$. Thus $v\in\ell^1(\Z)$, and we may normalize $v$
to satisfy $\langle v,v\rangle=1$ where $\langle\cdot,\cdot\rangle$ denotes
the duality (dot product) between $\ell^1(\Z)$ and $\ell^\infty(\Z)$.
The operator $L$ is Fredholm with index zero, with range
$$
\ran(L)=\{w\in\ell^\infty(\Z)\;|\;\langle v,w\rangle=0\},
$$
and its spectrum satisfies
$$
\sigma(L)\cap(0,\infty)=\emptyset.
$$
\exx
\vv

Note that the exponential estimate \eqref{expest} on $v_n$ implies
the absolute convergence of the sum \eqref{bsum} defining the quantity
$B$ in Condition B\@.

The proof of Proposition 1.4 will be given in Section $3$.
We remark that even if Condition A does not hold,
so $p$ in Proposition 1.4 is not unique,
we believe the statement of this result to be true, with $v$ depending on $p$.

To prove crystallographic pinning as in Theorem 1.1 we must prove
that the inequality \eqref{CP} is strict.
In our proof of Theorem 1.1 we shall assume that \eqref{CP}
is an equality for $\theta_0=0$ and that Condition A holds, and proceed to show that $B = 0$.
An important step in this proof is the analysis of the second order difference equation
\begin{equation}
X_{m+1} + X_{m-1} - 2X_m + (S+S^{-1}-2I)X_m - f(X_m,a_+(0))=0,\qquad m \in \Z, \label{diff}
\end{equation}
where $X_m\in\ell^\infty(\Z)$ for each $m$.
We write the vector $X_m$ in coordinate form as $X_m=\{x_{n,m}\}_{n\in\Z}$
where each $x_{n,m}$ is a scalar, and the operator $S$ acts on the
index $n$ so that $(SX_m)_n=x_{n+1,m}$.
In this way, equation \eqref{diff} is effectively a difference
equation in $x_{n,m}$ involving both indices $n$ and $m$. The function
$f(\cdot,a):\ell^\infty(\Z)\to\ell^\infty(\Z)$ in \eqref{diff}
is again a slight abuse of notation, where we evaluate the nonlinear scalar function
$f(\cdot,a):\R\to\R$ coordinatewise, that is,
$f(X_m,a)_n=f(x_{n,m},a)$.


Equation \eqref{diff} arises when
we consider the monotone traveling wave $\phi^\eps$ which propagates in
the direction $(1,\eps)$ with speed $c$, and we adjust the detuning parameter $a=a(\eps)$
so that $c$ is appropriately small.
This provides an infinite family of transition layers appearing at the integers $\xi = n$
(after an appropriate shift),
and to capture these layers we let $\phi^\eps_n(\zeta) = \phi^\eps(\eps\zeta +n)$.
Under this rescaling the differential-difference equation \eqref{2DWP} becomes the infinite system
\begin{equation}
-\frac{c}{\eps} ({\phi^\eps_n})'(\zeta) 
= \phi^\eps_n(\zeta+1) + \phi^\eps_n(\zeta-1) + \phi^\eps_{n+1}(\zeta) + \phi^\eps_{n-1}(\zeta)
- 4\phi^\eps_n(\zeta) - f(\phi^\eps_n(\zeta),a),
\label{layers}
\end{equation}
where each $\phi^\eps_n(\zeta)$ corresponds to the layer at $\xi=n$.
By choosing $a = a(\eps)$ so that $c = \eps^2$, the system \eqref{layers}
again develops an infinite family of transition layers spaced (in the limit) a unit distance apart,
that is, we have ``layers within layers.'' (This secondary scaling in our analysis takes the
form \eqref{twoscale}. Let us also remark
that solutions of PDE's with nested families of transition layers has been related
to solutions of variational problems on tori \cite{Bangert,Moser,Rabinowitz}.)
The vector $X_m = \{{\ds{\lim_{\eps \to 0}}} \phi_n^\eps(m)\}_{n \in \Z} \in \ell^\infty(\Z)$
captures the values of the plateaus between these secondary layers, and
if it is the case that \eqref{CP} is an equality, then the sequence $\{X_m\}_{m \in \Z}$
satisfies equation \eqref{diff}.
The vectors $X_m$ form a monotone heteroclinic orbit, having limits $X_{\pm \infty}$ which are approached
from above as $m\to-\infty$ and below as $m\to\infty$,
and in the event that Condition A holds, it is the case that
$X_{\infty} = SX_{-\infty}=p$, with $p$ as in Condition A\@. Furthermore, the limiting point $p$
possesses a two-dimensional center manifold as an equilibrium of the difference equation \eqref{diff}.
Theorem 1.1 is proved by showing that if
$B \ne 0$, then the dynamics on the center manifold has a Takens-Bogdanov normal form,
specifically \eqref{red2}, \eqref{gamma2}, albeit for a difference operator.
This in turn is shown to preclude the possibility of monotone limits to $p$ from both above and below,
and thus to preclude the monotone heteroclinic orbit $\{X_m\}_{m\in\Z}$. It follows that
$B=0$, as desired.



This paper is organized as follows. In Section 2 we construct the
heteroclinic solution $\{X_m\}_{m\in\Z}$ to equation \eqref{diff},
assuming that \eqref{CP} is an equality and that Condition A holds. In Section 3 we prove Proposition 1.4, which
develops information about the linearization of equation \eqref{diff} at the equilibrium $p$. This proof relies
on properties of resolvent positive operators. In Section 4 the center manifold reduction
along with associated shadowing properties is given for equation \eqref{diff}. In Section 5
we prove Theorem 1.1 by showing that $B=0$ must hold if \eqref{CP} is an equality.
Finally, in Section 6 we use transversality methods to prove Theorem 1.2,
namely that Conditions A and B hold generically for normal families of operators.


\section{Doubly Heteroclinic Orbits}

Throughout this section we assume that $f$ is a normal family, in particular
satisfying \eqref{bistable} and \eqref{normal}. We study the LDE \eqref{LDE},
keeping the notation as before.
\vv

\noindent {\bf Lemma 2.1.} \bxx
The inequalities  $-1<a_-(0)$ and $a_+(0)<1$ are strict.  That is, there is a $\delta > 0$ such that $c(a,(1,0))<0$ for $a\in(-1,-1+\delta]$
and $c(a,(1,0))>0$ for $a\in [1-\delta,1)$.
\exx
\vv

\noindent {\bf Proof.}
Without loss we show that $a_+(0)<1$. The conclusions about $c$ follow directly from this.  We proceed by contradiction.
Assume that $a_+(0)=1$.  Thus $c(a)\le 0$ for every $a\in(-1,1)$, where we denote
$c(a)=c(a,(\kappa,\sigma))$ with the horizontal direction $(\kappa,\sigma)=(1,0)$. Thus the
corresponding wave profile $\phi(\xi)=\phi(\xi;a,(1,0))$ satisfies
\begin{equation}
\phi(\xi+1)+\phi(\xi-1)-2\phi(\xi)-f(\phi(\xi),a)\ge 0
\label{L1}
\end{equation}
by \eqref{2DWP}. Without loss we may assume that
\begin{equation}
\phi(\xi)\le 0\hbox{ for }\xi<0,\qquad
\phi(\xi)\ge 0\hbox{ for }\xi>0,
\label{L2}
\end{equation}
after an appropriate translation of $\xi$. Now take any sequence $a_n\to 1$ and
let $\phi_n(\xi)$ be the corresponding wave profile, which is monotone in $\xi$.
By Helly's theorem there exists a subsequence
(which we still denote by $\phi_n(\xi)$) such that
\begin{equation}
\lim_{n\to\infty}\phi_n(\xi)=\phi_*(\xi)
\label{helly}
\end{equation}
exists for all but countably many $\xi\in\R$. Fixing such a subsequence, let
$$
S=\{\xi\in\R\;|\;\hbox{the limit (\ref{helly}) exists}\},\qquad
T=\{\xi\in\R\;|\;\xi+n\in S\hbox{ for every }n\in\Z\}.
$$
Clearly the complement $\R\setminus T$ of $T$ is countable.
Then from \eqref{L1} we see that
$$
\phi_n(\xi+1)+\phi_n(\xi-1)-2\phi_n(\xi)\ge\alpha_n=\min_{u\in[-1,1]}f(u,a_n),
$$
for every $\xi\in\R$, where the above equality serves as the definition of $\alpha_n$.
Clearly $\alpha_n\to 0$ as $n\to\infty$,
so upon taking this limit we have
\begin{equation}
\phi_*(\xi+1)+\phi_*(\xi-1)-2\phi_*(\xi)\ge 0,\qquad
\xi\phi_*(\xi)\ge 0,
\label{L3}
\end{equation}
for $\xi\in T$, where \eqref{L2} has been used for the second inequality. Morever, the first inequality
in \eqref{L3} is strict for those $\xi\in T$ for which $\phi_*(\xi)\ne\pm 1$. One sees this from \eqref{L1},
since $f(\phi_n(\xi),a_n)$ is strictly positive and bounded away from zero for large $n$.
Now take any $\xi_0\in T$ and denote
$\delta_n=\phi_*(\xi_0+n)-\phi_*(\xi_0+n-1)$. Then the first inequality in \eqref{L3} implies that
$\delta_{n+1}\ge\delta_n$ for every $n\in\Z$. Thus if $\delta_1>0$ we have
$$\phi_*(\xi_0+n)-\phi_*(\xi_0)=\sum_{k=1}^n\delta_k\ge n\delta_1\to\infty$$
as $n\to\infty$, which contradicts the boundedness of $\phi_*(\xi)$.
It follows that $\delta_1=0$, and thus $\phi_*(\xi+1)=\phi_*(\xi)$ for every $\xi\in T$.
But if $\xi\in T$ then also $\xi+n\in T$ for every $n$,
and so in fact $\phi_*(\xi+n)=\phi_*(\xi)$ holds, for every $n\in\Z$ and $\xi\in T$.
This, and the monotonicity of $\phi_*$ on $T$,
implies that $\phi_*(\xi)$ is constant on $T$. With this, the second inequality in \eqref{L3}
implies that $\phi_*(\xi)=0$ on $T$, and so the first inequality in \eqref{L3}
is an equality for $\xi\in T$. But this contradicts our assertion above that the first
inequality in \eqref{L3} is strict whenever $\phi_*(\xi)\ne\pm 1$.
\qed
\vv

\noindent {\bf Lemma 2.2.} \bxx
Let $q\in\ell^\infty(\Z)$ satisfy
\be
q_{n+1}+q_{n-1}-2q_n=f(q_n,a),\qquad n\in\Z,
\label{L4}
\ee
for some $a\in(-1,1)$. Assume also that $q_n\le q_{n+1}$ for every
$n\in\Z$. Then either
\be 
\lim_{n\to\pm\infty}q_n=\pm 1,
\label{qlim}
\ee
or else $q$ is a constant sequence with $q_n=q_0\in\{-1,a,1\}$ for every $n$.
\exx
\vv

\noindent {\bf Proof.}
Denoting $q_{\pm\infty}={\ds{\lim_{n\to\pm\infty}}}q_n$, which are finite
quantities, we have $f(q_{\pm\infty},a)=0$ from \eqref{L4}, and hence
$q_{\pm\infty}\in\{-1,a,1\}$. It is thus enough to prove that (a)
it is impossible for both $q_{-\infty}=-1$ and $q_\infty=a$ to hold simultaneously; and
(b) it is impossible for both $q_{-\infty}=a$ and $q_\infty=1$ to hold simultaneously.
Without loss we only prove (a), so assume to the contrary that $q_{-\infty}=-1$
and $q_\infty=a$. Noting that $f(q_n,a)\ge 0$ for every $n\in\Z$
and denoting $\delta_n=q_n-q_{n-1}$, we have from \eqref{L4}
that
\be
\delta_{n+1}=\delta_n+f(q_n,a)\ge\delta_n\ge 0.
\label{delta}
\ee
As $\{q_n\}_{n\in\Z}$ is not a constant sequence,
necessarily $\delta_m>0$ for some $m$. But then $\delta_{m+k}\ge\delta_m$
for every $k\ge 0$ by \eqref{delta}, hence
$$q_{m+k}-q_{m-1}=\sum_{i=0}^k\delta_{m+i}\ge (k+1)\delta_m\to\infty$$
as $k\to\infty$, contradicting boundedness of $q_{m+k}$.
\qed
\vv

\noindent {\bf Proposition 2.3.} \bxx
Assume that Condition A holds, with $p$ as stated there.
Then the $\ell^\infty(\Z)$-valued ODE
\be
-\Phi'=(S+S^{-1}-2I)\Phi-f(\Phi,a_+(0))
\label{2a}
\ee
admits a heteroclinic solution $\Phi:\R\to\ell^\infty(\Z)$ with limits
$$
\lim_{\zeta\to-\infty}\Phi(\zeta)=S^{-1}p,\qquad
\lim_{\zeta\to\infty}\Phi(\zeta)=p.
$$
Moreover, there exists such a solution
$\Phi(\zeta)=\{\phi_n(\zeta)\}_{n\in\Z}$ satisfying the monotonicity condition
$$
\phi_n(\zeta)<\phi_n(\zeta+\delta)\hbox{ for }\delta>0,
$$
for every $n\in\Z$ and $\zeta\in\R$.
\exx
\vv

\noindent {\bf Proof.}
Choose $b \in (-1,1) \setminus \{a\}$ so that also $b\ne p_j$ for every $j\in\Z$,
with $p$ as in Condition A\@.
We shall consider the monotone traveling wave $\phi^\eps = \phi^\eps(\xi) = \phi^\eps(\xi;a^\eps,(1,0))$ traveling in the direction $(1,0)$ with appropriately chosen $a^\eps$ and normalized
by a translation so that $\phi^\eps(0) = b$.  For sufficiently small $\eps$ choose $a^\eps$ so that $c = c(a^\eps,(1,0)) = \eps$.  Such a choice is  possible by the continuity of the function $c$ and by Lemma 2.1.  Note that ${\ds{\lim_{\eps \to 0}}} a^\eps = a_+(0)$.  

Let $\phi_n^\eps(\zeta) = \phi^\eps(\eps \zeta + n)$, and substitute into \eqref{WP} to obtain
\[ -(\phi_n^\eps)'(\zeta) = \phi_{n+1}^\eps(\zeta) + \phi_{n-1}^\eps(\zeta) - 2\phi^\eps_n(\zeta) - f(\phi_n^\eps(\zeta),a^\eps) \]
For each $n$, as $\eps$ varies $\phi_n^\eps$ constitute a bounded equicontinuous family and thus by Ascoli's theorem and a diagonalization argument, we have the limit
\[\lim_{\eps \to 0} \phi_n^\eps(\zeta) = \phi_n(\zeta)\] for some subsequence of $\eps$, holding uniformly on compact $\zeta$-intervals for every $n \in \Z$.  Letting $\Phi(\zeta) := \{\phi_n(\zeta)\}_{n \in \Z}$ we obtain \eqref{2a}.  
The monotonicity of the original wave profile $\phi^\eps$ together with the scaling $\zeta = \eps \xi + n$ implies the monotonicity condition
\begin{equation} \phi_n(\zeta) \le \phi_n(\zeta + \delta) \le \phi_{n+1}(\zeta) \qquad \delta > 0
\label{eq:mono1}
\end{equation}
holds for all $\zeta \in \R$ and all $n \in \Z$.  We now show that this monotonicity is strict, as in the statement
of the proposition.
Suppose to the contrary that $\phi_n(\zeta_0) = \phi_n(\zeta_0+\delta)$ for some $n\in\Z$ and some $\delta > 0$, and some $\zeta_0 \in \R$.
Since each $\phi_n$ is non-decreasing we have $\phi_n(\zeta) = \phi_n(\zeta_0)$ for
every $\zeta \in [\zeta_0,\zeta_0+\delta]$.  It is then a consequence of \eqref{2a} that
$\phi_{n+1}(\zeta) + \phi_{n-1}(\zeta)$ is constant in $[\zeta_0,\zeta_0 + \delta]$, and so
both $\phi_{n\pm 1}(\zeta)$ are constant in that interval as these functions are nondecreasing.
By induction we see that the vector $\Phi(\zeta)$ is in fact constant in this interval,
and so \eqref{2a} at $\zeta=\zeta_0$
rewrites as $(S+S^{-1}-2I)\Phi(\zeta_0)-f(\Phi(\zeta_0),a_+(0)) = 0$. Thus the point
$\Phi(\zeta_0)\in\ell^\infty(\Z)$ is an equilibrium of equation \eqref{2a}, and by uniqueness
we must have $\Phi(\zeta)=\Phi(\zeta_0)$ for every $\zeta\in\R$.
Lemma 2.2 now applies to the vector $q=\Phi(\zeta_0)$
due to the monotonicity of $\phi_n(\zeta_0)$ in $n$. Because
$q_0=\phi_0(0)=b\not\in\{-1,a,1\}$ it follows that \eqref{qlim} holds.
Thus by the uniqueness (up to translation) of the vector $p$ in Condition~A,
we have that $q=S^jp$ for some $j\in\Z$. But then $b=q_0=p_j$ contradicts the
choice of $b$. Thus the first inequality
in  \eqref{eq:mono1} is strict, as desired. Note that the second inequality in \eqref{eq:mono1} is also strict.

We now establish the boundary conditions for $\Phi(\zeta)$.
It is a consequence of \eqref{eq:mono1} that the limits $\phi_n(\pm \infty) = {\ds{\lim_{\zeta \to \pm \infty}}} \phi_n(\zeta)$ exist.
Taking these limits in \eqref{2a} gives
$(S + S^{-1} - 2I)\Phi(\pm\infty) - f(\Phi(\pm\infty),a_+(0)) = 0$.
Again Lemma 2.2 applies, and this time due to the strict monotonicity
we have \eqref{qlim} for both vectors $\Phi(\pm\infty)$.
As $\Phi(-\infty)\le \Phi(\infty)$ and $\Phi(-\infty)\ne\Phi(\infty)$ it follows that
$\Phi(-\infty) = S^jp$ and $\Phi(\infty)=S^k p$ for some integers $j<k$, with $p$ given by Condition A\@.
However, we also have that $\Phi(\infty)\le S\Phi(-\infty)$
from the second inequality in \eqref{eq:mono1}, and so $k\le j + 1$. Thus $k = j+1$.
By reindexing $n$, namely by replacing $\Phi(\zeta)$ with $S^{-k}\Phi(\zeta)$, we obtain the desired boundary conditions.
\qed
\vv

\noindent {\bf Proposition 2.4.} \bxx
Assume that the inequality \eqref{CP} is an equality at $\theta_0=0$.
Also assume that Condition A holds, with $p$ as stated there.
Then the $\ell^\infty(\Z)$-valued difference equation \eqref{diff} admits a heteroclinic orbit
$\{X_m\}_{m\in\Z}$ with limits
\be
\lim_{m\to\pm\infty}X_m=X_{\pm\infty},\qquad
X_{-\infty}=S^{-1}p,\qquad
X_\infty=p,
\label{xlims}
\ee
in the space $\ell^\infty(\Z)$.
Moreover $X_m=\{x_{n,m}\}_{n\in\Z}$ satisfies the lexicographic monotonicity property
\begin{equation} \label{2het}
p_{n-1}=x_{n,-\infty}<x_{n,m}<x_{n,m+1}<x_{n,\infty}=p_n,\qquad n,m\in\Z.
\end{equation}
\exx
\vv

\noindent {\bf Proof.}
We shall consider the monotone traveling wave
$\phi^\eps=\phi^\eps(\xi)=\phi(\xi;a^\eps,(1,\eps))$,
which propagates in direction $(\kappa,\sigma)=(1,\eps)$, taking an
appropriately chosen $a=a^\eps$.
Note that the angle $\theta$ corresponding to this direction is $\theta=\arctan\eps$,
and so the upper value of the pinning range is located at $a=a_+(\theta)=a_+(\arctan\eps)$.

To be precise, fix a sequence $\eps_k\to 0$ with $\eps_k\ne 0$ such that
$$
\lim_{k\to\infty}a_+(\theta_k)=a_+(0)
$$
for $\theta_k=\arctan\eps_k$, by the equality assumption in \eqref{CP}. For the
remainder of this proof we shall select $\eps$ only from the sequence $\eps_k$,
suppressing the index $k$ for notational simplicity.
For sufficiently small $\eps$ choose
$a = a^\eps$ so that the wave speed is given by $c(a^\eps,(1,\eps)) = \eps^2$.
Such a choice is possible by the continuity of the function $c=c(a,(\kappa,\sigma))$ and
by Lemma 2.1. Indeed, one has $c(a_+(\arctan\eps),(1,\eps))=0$ while
$c(1-\delta,(1,\eps))>0$ is bounded away from zero as $\eps\to 0$, by Lemma 2.1.
Let $\bar{a} = {\ds{\lim_{\eps \to 0}}} a^\eps$, where the limit is taken along a subsequence if necessary.
Certainly $a^\eps > a_+(\arctan\eps)$ for all $\eps$
and therefore
$\bar{a} \ge {\ds{\lim_{\theta \to 0}}} a_+(\arctan\eps)=a_+(0)$.
Similarly, continuity of the function $c$ implies that $c(\bar{a},(1,0)) = 0$,
so $\bar{a} \le a_+(0)$. Thus $\bar{a} = a_+(0)$.

Now fix $b\in (p_{-1},p_0)$ so that $b\ne a_+(0)$. Of course $b\ne p_n$ for every $n$.
By translating the argument
$\xi$ of $\phi^\eps$, we may assume without loss that $\phi^\eps(0) = b$.
Let
\be
\phi^\eps_{n,m}(\zeta) = \phi^\eps(\eps^2 \zeta + \eps m + n).
\label{twoscale}
\ee
As $\phi^\eps(\xi)$
satisfies equation \eqref{2DWP} with $c=\eps^2$,
with $a=a^\eps$, and with $(\kappa,\sigma)=(1,\eps)$, we see that
the functions $\phi_{n,m}^\eps(\zeta)$ satisfy the LDE
\begin{equation}
-(\phi_{n,m}^\eps)'
= \phi^\eps_{n,m+1}+\phi^\eps_{n,m-1}+\phi^\eps_{n+1,m}+\phi^\eps_{n-1,m}
-4\phi^\eps_{n,m} - f(\phi^\eps_{n,m},a(\eps)).
\label{nmLDE}
\end{equation}
Note in particular that all arguments of the functions $\phi^\eps_{i,j}$ in \eqref{nmLDE}
are evaluated at the same point $\zeta$, and so \eqref{nmLDE} is an infinite-dimensional ODE\@.
By Ascoli's theorem and a diagonalization argument, we have the limit
$$\lim_{\eps\to 0}\phi^\eps_{n,m}(\zeta)=\phi_{n,m}(\zeta)$$
for some subsequence of $\eps$, holding uniformly on compact $\zeta$-intervals
for every $n$ and $m$. 
The limiting function $\phi_{n,m}$ satisfies the same LDE \eqref{nmLDE}, but with $a_+(0)$
in the argument of $f$. Of course, $\phi_{0,0}(0)=b$.
The monotonicity of the original wave profile $\phi^\eps(\xi)$,
and the particular scaling $\xi=\eps^2\zeta+\eps m+n$ of the argument, imply that
our limiting functions enjoy the lexicographic montonicity condition
\begin{equation}
\phi_{n,m}(\zeta) \le \phi_{n,m}(\zeta+\delta) \le \phi_{n,m+k}(\zeta) \le \phi_{n+1,m}(\zeta-\delta)\hbox{ for }
k\ge 1\hbox{ and }\delta>0, \label{philex}
\end{equation}
and every $n,m\in\Z$ and $\zeta\in\R$.
Indeed, one sees that for any given such $k$, $\delta$, $n$, $m$, and $\zeta$,
the inequalities \eqref{philex} hold for $\phi^\eps_{n,m}$ for small $\eps$.
Define now
$$x_{n,m}=\lim_{\zeta\to\infty}\phi_{n,m}(\zeta),\qquad
X_m=\{x_{n,m}\}_{n\in\Z}\in\ell^\infty(\Z).$$
Certainly these limits exist, and of course $|x_{n,m}|\le 1$.
Additionally ${\ds{\lim_{\zeta\to\infty}}}\phi'_{n,m}(\zeta)=0$
from the differential equation \eqref{nmLDE}. From this one easily sees that taking the limit $\zeta\to\infty$
in \eqref{nmLDE} yields the
difference equation \eqref{diff}, which may also be written as
\be
(\Delta x)_{n,m}=f(x_{n,m},a_+(0)),\qquad
(n,m)\in\Z^2,
\label{xlap}
\ee
where $\Delta$ is the discrete laplacian \eqref{2dlap} on $\Z^2$.
The ordering $x_{n,m}\le x_{n,m+k}\le x_{n+1,m}$ for $k\ge 1$,
which follows from \eqref{philex}, implies that
$x_{n,\pm\infty}={\ds{\lim_{m\to\pm\infty}}}x_{n,m}$ exist and that
\be
x_{n,-\infty}\le x_{n,m}\le x_{n,m+1}\le x_{n,\infty},
\qquad
x_{n,\infty}\le x_{n+1,-\infty},
\label{xlex}
\ee
for every $n$ and $m$. One sees that the vectors
$X_{-\infty}=\{x_{n,-\infty}\}_{n\in\Z}$ and $X_{\infty}=\{x_{n,\infty}\}_{n\in\Z}$
both satisfy the conditions of Lemma 2.2 with $a=a_+(0)$.

Also observe that
\be
x_{0,-1}\le b\le x_{0,0}.
\label{bvalue}
\ee
The second inequality in \eqref{bvalue} holds because
$\phi_{0,0}(\zeta)\ge\phi_{0,0}(0)=b$ for $\zeta\ge 0$.
The first inequality in \eqref{bvalue} follows by taking
$(n,m)=(0,-1)$ with $k=1$, and $\zeta=0$, in the second inequality
of \eqref{philex}, to give $\phi_{0,-1}(\delta)\le\phi_{0,0}(0)=b$.
One then lets $\delta\to\infty$.

We next show that the first three inequalities in \eqref{xlex} are strict.
Suppose that $x_{n,m}=x_{n,m+1}$ for some $n$ and $m$. Then by equation \eqref{xlap}
$$\ba{lcl}
0 &\bb = &\bb (\Delta x)_{n,m+1}-(\Delta x)_{n,m}\\
\\
&\bb = &\bb (x_{n+1,m+1}-x_{n+1,m})
+(x_{n-1,m+1}-x_{n-1,m})
+(x_{n,m+2}-x_{n,m-1}).
\ea$$
Each of the three differences in the final line of the above equation
are nonnegative by \eqref{xlex}, and so we have that
$$
\begin{array}{l}
x_{n+1,m}=x_{n+1,m+1},\qquad
x_{n-1,m}=x_{n-1,m+1},\\
\\
x_{n,m-1}=x_{n,m}=x_{n,m+1}=x_{n,m+2}.
\end{array}
$$
It is clear by repeating the procedure that we may conclude that
\be
x_{n,m}=x_{n,0}
\label{mindep}
\ee
for every $n$ and $m$, namely, that $X_m=X_0$
is constant in $m$. Thus $X_{\pm\infty}=X_0$, and as noted this
vector satisfies the conditions of Lemma 2.2 with $a=a_+(0)$.
We also have, from \eqref{bvalue} and \eqref{mindep}, that
$x_{0,0}=b\not\in(-1,a_+(0))\cup(a_+(0),1)$,
and so it follows from Lemma 2.2 that
${\ds{\lim_{n\to\pm\infty}}}x_{n,0}=\pm 1$. By the uniqueness assumption on $p$
in the statement of the proposition, we have that $X_0=S^kp$ for some $k\in\Z$.
But then $b=x_{0,0}=p_k$, which contradicts the choice of $b$. We thus conclude
that $x_{n,m}<x_{n,m+1}$ for every $n$ and $m$, and hence the first three inequalities
in \eqref{xlex} are strict.

One thus has
$x_{n,-\infty}<x_{n,\infty}\le x_{n+1,-\infty}<x_{n+1,\infty}$ for every $n$,
and this implies that $x_{n,-\infty}<x_{n+1,-\infty}$ and
$x_{n,\infty}<x_{n+1,\infty}$.
Applying Lemma 2.2 again, to both $X_{-\infty}$ and $X_\infty$,
we conclude that each of these
vectors are shifts of $p$, say
$$X_{-\infty}=S^jp,\qquad X_\infty=S^kp,$$
for some $j$ and $k$. As $x_{n,-\infty}<x_{n,\infty}$ we have that
$S^jp\le S^kp$ with $S^jp\ne S^kp$, and so $j<k$.
As $x_{n,\infty}\le x_{n+1,-\infty}$ we have that $S^kp\le S^{j+1}p$, hence
$k\le j+1$. Thus $j=k-1$. Further,
$$p_{k-1}=(S^{k-1}p)_0=(S^jp)_0=x_{0,-\infty}<x_{0,0}=b<x_{0,\infty}=(S^kp)_0=p_k,$$
and so $k=0$ by the choice of $b$. This implies that $x_{n,-\infty}=p_{n-1}$
and $x_{n,\infty}=p_n$, to give the result.
\qed

\section{The Proof of Proposition 1.4}

For each pair $(q,a) \in \ell^\infty(\Z) \times (-1,1)$ we may linearize equation \eqref{L4}
about $q$ to obtain an associated operator
$L_q(a) \in \mathcal{L}(\ell^\infty(\Z))$ given by
\be L_q(a) = S+S^{-1} - 2I - f'(q,a), \label{Ldef}
\ee
where $f'(q,a)$ denotes the diagonal operator with entries $f'(q_n,a)$.
Proposition 1.4 concerns the operator $L_q(a)$ with
$a = a_+(0)$ and with $q=p$ as in \eqref{pdiff}, \eqref{BM}.
We first examine the spatial dynamics associated with the difference equation $L_q(a)v = 0$.

Generally, we shall denote $L_q=L_q(a)$ when the value of $a$ is clear.
\vv

\noindent {\bf Lemma 3.1.} \bxx
Let $D\in\mathcal{L}(\ell^\infty(\Z))$ be the diagonal operator
$D=\diag\{d_n\}_{n\in\Z}$, where we assume the existence of the limits and the sign conditions
$$\lim_{n\to\pm\infty} d_n=d_{\pm\infty},\qquad d_{\pm\infty}>0.$$
Then the operator $T=S+S^{-1}-2I-D$ is Fredholm of index zero.
Moreover, the dimension of the kernel of $T$ is either zero or one.
Furthermore, elements $v \in \ker(T)$ of the kernel
enjoy the estimate
\be
|v_n|\le K\mu^{|n|},\qquad n\in\Z,
\label{expest2}
\ee
for some $K>0$ and $0<\mu<1$, and in addition,
\be
v_n\ne 0\hbox{ for all large }|n|
\label{vnz}
\ee
for nontrivial elements of the kernel. If $v$ is a nontrivial element of the kernel then
\be
\ran(T)=\{w\in\ell^\infty(\Z)\;|\;\langle v,w\rangle=0\},
\label{kerran2}
\ee
with $\ran(T)=\ell^\infty(\Z)$ if $\ker(T)=\{0\}$.

If $q \in \ell^\infty(\Z)$ satisfies ${\ds{\lim_{n\to\pm\infty}}}q_n=\pm 1$
and $a \in (-1,1)$, then all the above claims hold for the operator $T=\lambda I-L_q$, where $L_q=L_q(a)$, provided that
$\lambda > \max\{-f'(1,a),-f'(-1,a)\}$.

Finally, the same conclusions hold if we consider the operator $T$, or $\lambda I-L_q$, lying in the space $\LL(Q_0)$
instead of $\LL(\ell^\infty(\Z))$, where
$$
Q_0=\{x\in\ell^\infty(\Z)\;|\;\lim_{n\to\pm\infty}x_n=0\}.
$$
\exx
\vv

\noindent {\bf Proof.}
The arguments below apply equally to operators in the space $\LL(\ell^\infty(\Z))$,
and in the space $\LL(Q_0)$.

Consider the second order difference equation
\be
x_{n+1}+x_{n-1}-(2+d_n)x_n=0
\label{diffeq}
\ee
associated to the operator $T$.  This difference equation
is asymptotically hyperbolic, that is, all roots $\mu$ of the limiting characteristic equations
\be
\mu+\mu^{-1}-(2+d_{\pm\infty})=0
\label{chareq}
\ee
satisfy $|\mu|\ne 1$. Indeed, one easily checks this by noting the left-hand side
of \eqref{chareq} is negative whenever $\mu=e^{i\theta}$.
In fact the two roots of this equation have the form $\mu=\mu_\pm$ and $\mu=\mu_\pm^{-1}$,
satisfying $0<\mu_\pm<1<\mu_\pm^{-1}$, for each choice of sign $+$ or $-$, and
so this equation is of saddle type.

Asymptotic hyperbolicity of equation \eqref{diffeq} now implies that this equation admits dichotomies whose stable and unstable subspaces are each
one-dimensional, both for $n\to\infty$
and for $n\to-\infty$.
This and the above remarks about the roots of \eqref{chareq}
in turn imply that the kernel of $T$ is at most one-dimensional, that any kernel element necessarily decays at the rate
\eqref{expest2} for some $K>0$ and $0<\mu<1$, and that $T$ is Fredholm.

In particular, the proof that $T$ is Fredholm follows by adapting the techniques of Palmer \cite{Palmer} for ODE's to the setting of difference equations.
Palmer proved that a differential operator of the form
$Lx(t)=\dot x(t)-A(t)x(t)$, where $A(t)$ is a bounded continuous coefficient matrix, is Fredholm if and only if it admits exponential dichotomies on both half lines.
(Here the operator $L$ acts from the space of bounded $C^1$ functions on the line with bounded
first derivative, into the space of bounded continuous functions.)
The techniques of Palmer's proof carry over to difference operators of the form \eqref{diffeq},
with appropriate modifications. We omit the details.

To prove that the Fredholm index of $T$ is zero, consider the operator $T-\lambda I$. For
every $\lambda\ge 0$ this operator satisfies the same conditions in the statement of the
lemma as $T$ does, and so it too is Fredholm, necessarily with the same index as $T$.
But $T-\lambda I$ is invertible for $\lambda>\|T\|$,
and so has index zero. Thus the index of $T$ also is zero. In particular, if $\ker(T)=\{0\}$
then $\ran(T)=\ell^\infty(\Z)$.

To prove \eqref{kerran2} in the case that $\ker(T)$ contains a nontrivial
element $v$, note that if $w\in\ran(T)$, say $w=Ty$, then
$\langle v,w\rangle=\langle v,Ty\rangle
=\langle T^*v,y\rangle=0$,
where $T^*\in\mathcal{L}(\ell^1(\Z))$ has the same formula as does $T$,
but considered in the space $\ell^1(\Z)$.
This now proves an inclusion for the formula \eqref{kerran2}.
But both spaces in this formula have the same codimension, in particular
$\codim\ran(T)=\dim\ker(T)=1$,
as the Fredholm index of $T$ is zero. We conclude that
\eqref{kerran2} is valid, as claimed.

To prove that \eqref{vnz} holds for any nontrivial element of $\ker(T)$, consider
such an element $v$. Suppose for some $n\ge n_0$ that $v_{n+1}>v_n\ge 0$,
where $n_0$ is large enough that $d_m\ge0$ for every $m>n_0$. Then from the difference equation \eqref{diffeq}
at $n+1$ we have that $v_{n+2}-v_{n+1}\ge v_{n+1}-v_n$, and so $v_{n+2}>v_{n+1}\ge0$. Upon repeating this procedure, we
conclude that $v_{n+k}-v_{n+k-1}$ is a nondecreasing and strictly positive sequence for $k\ge 1$.
But this forces $v_m\to\infty$ as $m\to\infty$, which is false. Thus the assumption
that $v_{n+1}>v_n\ge 0$ for some $n\ge n_0$ is impossible. Similarly, $v_{n+1}<v_n\le 0$ cannot occur for any $n\ge n_0$.
Thus if $v_n=0$ for some $n\ge n_0$, necessarily $v_{n+1}=0$. However, in light of the difference equation \eqref{diffeq},
this forces $v_m=0$ for every $m\in\Z$, which contradicts the assumption that $v$ is nontrivial.
We thus conclude that $v_n\ne 0$ for all
$n\ge n_0$. Similarly, $v_n\ne 0$ for all sufficiently large negative $n$, as claimed.

To establish the penultimate claim of the lemma, we note that
the operator $L_q-\lambda I$ has the same form as $T$ in the statement of the lemma, with
$$d_{\pm\infty}=f'(q_{\pm\infty},a)+\lambda=f'(\pm 1,a)+\lambda.$$
Thus $L_q-\lambda I$ is Fredholm with index zero so long as $\lambda > \max\{-f'(1,a),-f'(-1,a)\}$, as desired.
\qed
\vv

The proof of Proposition 1.4 relies on the theory of resolvent positive operators, which we now briefly outline.

Recall that a closed convex subset $\mathcal{K}$ of a Banach Space $\X$
is called a {\bf cone} if it is closed under positive linear combinations
and if also $\mathcal{K} \cap (-\mathcal{K}) = \{0\}$.
A cone $\mathcal{K}$ is a called a {\bf total cone} if in addition, its linear span
$\{x_1-x_2\;|\;x_1,x_2\in\mathcal K\}$ is dense in $\X$.

Let $\mathcal K\subseteq\X$ be a cone. An element $x \in \X$ is called {\bf positive}
with respect to $\mathcal K$, denoted $x\ge 0$, if $x \in \mathcal{K}$,
and we write $x_1\ge x_2$ to mean $x_1-x_2\ge 0$, for $x_1,x_2\in\X$. The relation $\ge$
thus defines a partial order on $\X$.
A linear operator $T\in\mathcal{L}(\X)$ is called {\bf positive} if
$T(\mathcal{K}) \subseteq \mathcal{K}$, that is, $Tx\ge 0$ whenever $x\ge 0$.
We similarly write $T\ge 0$ to denote that $T$ is a positive operator,
and $T_1\ge T_2$ when $T_1-T_2\ge 0$.
An operator $T$ is called {\bf resolvent positive} if there exists $\Lambda\in\R$
such that $\lambda\in\rho(T)$ and
$(\lambda I - T)^{-1}\ge 0$ whenever $\lambda\ge\Lambda$.

Resolvent positive operators satisfy the following nice property \cite{JMP-CP,Nussbaum,Thieme}.
A proof of the following result can be found, for example, in Proposition 4.7 of \cite{JMP-CP}.
\vv

\noindent {\bf Proposition 3.2.} \bxx
Let $\mathcal{K}$ be a total cone in a Banach Space $\X$ and
let $T \in \mathcal{L}(\X)$ be resolvent positive with respect to $\mathcal{K}$.
Then
$$\sigma(T)\cap\R\ne\emptyset.$$
Further, let
$$\lambda_0 = \sup (\sigma(T) \cap \R).$$
Then either $\lambda_0$ is an eigenvalue of $T$ and there exists
an associated positive eigenvector $v$, that is
$$Tv=\lambda_0v,\qquad v\in\mathcal{K}\setminus\{0\},$$
or else the operator $\lambda_0 I-T$ is not Fredholm.
\exx
\vv

In what follows, whenever we take $\X=\ell^\infty(\Z)$, then we let
$\mathcal{K}$ be the cone of non-negative sequences,
namely those $x=\{x_n\}_{n\in\Z}$ with $x_n\ge 0$ for every $n$.
\vv

\noindent {\bf Lemma 3.3.} \bxx
If $T\in\mathcal{L}(\X)$ satisfies $T\ge -\sigma I$ for some $\sigma\in\R$,
then $T$ is resolvent positive. In particular, the operator $L_q(a)$
in equation \eqref{Ldef} is of this form for any
$q \in \ell^\infty(\Z)$ and $a \in (-1,1)$, and hence is resolvent positive.
\exx
\vv

\noindent {\bf Proof.}
For a general operator $T$ as in the statement of the lemma, write
$T=B-\sigma I$ where $B\ge 0$. Then
$$(\lambda I-T)^{-1}=((\lambda+\sigma)I-B)^{-1}
=\frac{1}{\lambda+\sigma}\sum_{k=0}^\infty \frac{B^k}{(\lambda+\sigma)^k}$$
converges for $\lambda+\sigma>\|B\|$, and is a sum of positive operators
hence is positive.

For the specific operator $L_q$ we see that taking
$\sigma \ge 2 + \max\{ f'(u,a) \; | \; -1 \le u \le 1\}$
gives $L_q\ge -\sigma I$, as desired.
\qed
\vv

\noindent {\bf Proposition 3.4.} \bxx
Suppose that $q \in \ell^\infty(\Z)$ and $a\in(-1,1)$
are as in the statement of Lemma 3.1.
Then $\sigma(L_q)\cap\R\ne\emptyset$. Assume additionally that the quantity
\be
\lambda_0 = \sup (\sigma(L_q) \cap \R)
\label{33a}
\ee
satisfies
\be\label{fredineq}
\lambda_0 > \max\{-f'(1,a),-f'(-1,a)\}.
\ee
Then $\lambda_0$ is an eigenvalue of $L_q$ and $\lambda_0I-L_q$ has a one-dimensional
kernel, and so
$$
\ker(\lambda_0I-L_q)=\{av\;|\;a\in\R\}
$$
for some $v\in\ell^\infty(\Z)\setminus\{0\}$.
Further, the eigenvector $v$ can be chosen to satisfy
$v_n > 0$ for all $n \in \Z$, and there are constants $K > 0$ and $0 < \mu < 1$
such that $v_n \le K\mu^{|n|}$ for all $n \in \Z$, and so also $v\in\ell^1(\Z)$.
Finally, the operator $\lambda_0 I-L_q$ is Fredholm with index zero, with range
\be
\ran(\lambda_0I-L_q)=\{w\in\ell^\infty(\Z)\;|\;\langle v,w\rangle=0\}.
\label{kerran}
\ee
\exx
\vv

\noindent {\bf Proof.}
The operator $L_q$ is resolvent positive by Lemma 3.3 and so
Proposition 3.2 implies that $\sigma(L_q)\cap\R\ne\emptyset$.
Now assume that the inequality \eqref{fredineq} holds.
Then $\lambda_0I-L_q$ is Fredholm with index zero, by Lemma 3.1.
Thus Proposition 3.2 implies that $\lambda_0$
is an eigenvalue of $L_q$ with an associated eigenvector $v\in\mathcal{K}\setminus\{0\}$,
that is, with $v\ne 0$, and $v_n\ge 0$ for every $n\in\Z$.
It follows from Lemma 3.1 that $\lambda_0$ is simple and that $v$ enjoys exponential decay in $n$,
and that moreover the formula \eqref{kerran} for the range holds.

It remains to prove that the strict inequality $v_n > 0$ holds for every $n$. If this is false, then there exists
some $k\in\Z$ for which $v_k=0$, and with either $v_{k+1}>0$ or $v_{k-1}>0$.
In any case, both $v_{k\pm 1}\ge 0$. But then
$$0=((\lambda_0I-L_q)v)_k=\lambda_0v_k-v_{k+1}-v_{k-1}+(2+f'(q_k,a))v_k<0,$$
a contradiction.
\qed
\vv

\noindent {\bf Proof of Proposition 1.4.}
It follows from Proposition 3.4 that it suffices to show that
$\lambda_0=0$ for the quantity in \eqref{33a}, taking $q=p$
and $a=a_+(0)$, and thus with the operator $L=L_p(a_+(0))$
as in equation \eqref{lop}.

Let us first show that $\lambda_0\le 0$. Assuming to the contrary that
$\lambda_0>0$, let $v\in\ell^\infty(\Z)$ be the associated positive
eigenvector as in Proposition 3.4.
Also let $\Psi(\zeta)=S\Phi(\zeta)-p$, with $\Phi(\zeta)$
as in Proposition 2.3, so that $\Psi(\zeta)\ge 0$ with
$\Psi(\zeta)\ne 0$. Then from equation \eqref{2a} we have
$$
-\Psi'=(S+S^{-1}-2I)\Psi-f(\Psi+p,a_+(0))+f(p,a_+(0))
$$
which can be rewritten as
$$
\begin{array}{lcl}
-\Psi' &\bb = &\bb L\Psi-E(\zeta),\\
\\
E(\zeta) &\bb = &\bb
f(\Psi(\zeta)+p,a_+(0))-f(p,a_+(0))-f'(p,a_+(0))\Psi(\zeta).
\end{array}
$$
Now $v\in\ell^1(\Z)$ by the exponential decay of $v_n$, and so we may take the
inner (duality) product $\langle\cdot,\cdot\rangle$ of the above equation with $v$. Denoting this by
$w(\zeta)=\langle v,\Psi(\zeta)\rangle$, we have that
\be
\begin{array}{lcl}
-w'(\zeta) &\bb = &\bb
\langle v,L\Psi(\zeta)\rangle-\langle v,E(\zeta)\rangle\\
\\
&\bb = &\bb
\langle L^*v,\Psi(\zeta)\rangle-\langle v,E(\zeta)\rangle
=\lambda_0 w(\zeta)-e(\zeta),
\end{array}
\label{3a}
\ee
where $e(\zeta)=\langle v,E(\zeta)\rangle$. Also note that $w(\zeta)>0$ for every $\zeta$. Now
$$
E_n(\zeta)=f(\psi_n(\zeta)+p_n,a_+(0))-f(p_n,a_+(0))-f'(p_n,a_+(0))\psi_n(\zeta)
$$
and so there exists $K>0$ such that
$$
|E_n(\zeta)|\le K\psi_n(\zeta)^2\le K\|\Psi(\zeta)\|\psi_n(\zeta)
$$
with $\|\cdot\|$ denoting the norm in $\ell^\infty(\Z)$. Thus
$$
|e(\zeta)|\le\sum_{n=-\infty}^\infty v_n|E_n(\zeta)|
\le K\|\Psi(\zeta)\|\sum_{n=-\infty}^\infty v_n\psi_n(\zeta)
=K\|\Psi(\zeta)\|w(\zeta),
$$
where the positivity of $v$ and of $\Psi(\zeta)$ is crucial in obtaining
this inequality. From this and from \eqref{3a} we have that
$$
-w'(\zeta)\ge (\lambda_0-K\|\Psi(\zeta)\|)w(\zeta),
$$
and in particular, $-w'(\zeta)\ge\frac{\lambda_0}{2}w(\zeta)$ for all
sufficiently negative $\zeta$, since $\Psi(\zeta)\to 0$ as $\zeta\to-\infty$.
However, this is impossible as $w(\zeta)\to 0$ as $\zeta\to-\infty$ with
$w(\zeta)>0$ for every $\zeta$. This contradiction proves that $\lambda_0\le 0$.

To complete the proof that $\lambda_0=0$, we prove that $\lambda_0\ge 0$.
Assume to the contrary that $\lambda_0<0$. Then $0\not\in\sigma(L)$ from
the definition of $\lambda_0$ and so $L$ is invertible. This fact together
with the implicit function theorem allows us to solve equation \eqref{L4}
for $q=q(a)\in\ell^\infty(\Z)$ depending smoothly on $a$, for $a$ near
$a_+(0)$, and with $q(a)=p$ at $a=a_+(0)$. To be specific, assume that
such $q(a)$ is defined for $a\in U$, where $U\subseteq(-1,1)$ is an open
interval containing the point $a_+(0)$.
(We reserve the right below to decrease $U$ by taking a smaller neighborhood
of $a_+(0)$.) Let us recall here the
strict inequalities $-1<a_+(0)<1$ from Lemma 2.1,
which will also be used below.

It is enough to show that $q_n(a)$ depends monotonically on $n$,
that is, $q_n(a)\le q_{n+1}(a)$ for every $n\in\Z$ and $a\in U$. Indeed,
if this holds then for every $a\in U$ we have a monotone solution
of equation \eqref{2DWP} at $(\kappa,\sigma)=(1,0)$
with $c=0$, and this implies that $U$ is a subset of the pinning region
$(a_-(0),a_+(0))$. However, this contradicts the fact that $U$ is a
neighborhood of $a_+(0)$.

Let us prove strict monotonicity in $n$ for $n\ge 0$,
the proof of monotonicity for $n\le 0$ being similar. Fix a sufficiently
large integer $N$ and sufficiently small $\eps>0$, such that
\be
p_N>\frac{1+a_+(0)}{2}>a+\eps,\qquad \|q(a)-p\|<\eps,\qquad p_{n+1}-p_n>2\eps,
\label{ineqs}
\ee
holds for every $a\in U$ and for $0\le n<N$. One easily does this by first choosing
$N$ to satisfy the first inequality in \eqref{ineqs}, then choosing $\eps$, and then
possibly decreasing the neighborhood $U$. We claim that $q_{n+1}(a)> q_n(a)$ for
every $n\ge 0$ and $a\in U$. For such $a$, first observe that
$$
q_{n+1}(a)-q_n(a)>p_{n+1}-p_n-2\eps>0,\qquad 0\le n<N,
$$
by \eqref{ineqs}. Also observe that
\be
q_n(a)>p_n-\eps\ge p_N-\eps>\frac{1+a_+(0)}{2}-\eps>a,\qquad n\ge N,
\label{largen}
\ee
again by \eqref{ineqs}.
Now assume for some $a\in U$ that there exists $m\ge N$ such that
$q_{m+1}(a)-q_m(a)\le 0$. Take $m$ to be the least such integer, and so
$q_m(a)-q_{m-1}(a)> 0$. Then from equation \eqref{L4} at $n=m$ we
have that $f(q_m(a),a)<0$, and so from \eqref{largen} we have that
$q_m(a)\in(a,1)$. Two cases now arise, and we shall rule out each of
them in turn. For the first case, there
exists $j>m$ such that $q_{j+1}(a)-q_j(a)>0$. Take $j$ to
be the least such integer, and so $q_j(a)-q_{j-1}(a)\le 0$, which
with equation \eqref{L4} implies that $f(q_j(a),a)>0$. However,
from the choice of $j$ and from \eqref{largen} it follows that
$a<q_j(a)\le q_m(a)<1$, which implies that $f(q_j(a),a)<0$, a contradiction.
We now consider the second case, in which $q_{j+1}(a)-q_j(a)\le 0$
for every $j>m$. Then the limit $q_j(a)\to q_\infty(a)$ exists
as a nonincreasing sequence for $j\ge m$, and upon taking the
limit $n\to\infty$ in equation \eqref{L4} we have that
$f(q_\infty(a),a)=0$, hence $q_\infty(a)\in\{-1,a,1\}$. But \eqref{largen}
implies that $q_\infty(a)>a$, and so $q_\infty(a)=1$. On the other hand,
$q_\infty(a)\le q_m(a)<1$, and this is a contradiction. With this, the proof of
the proposition is complete.
\qed
\vv

The fact that $\lambda_0 = 0$ will allow us to make a center manifold reduction.
The positivity of the kernel element $v$ will allow us to make use of the
center manifold reduction to induce a sign condition on $f$.  The next section is devoted to constructing the center manifold.

\section{A Center Manifold Reduction}

So long as Condition A holds and \eqref{CP} is an equality at $\theta_0=0$,
then Proposition 2.4 guarantees that the difference equation
\eqref{diff} admits a monotone increasing solution $X_m$ which converges to $p$ as $m \to \infty$.
Similarly, the solution of \eqref{diff} given by $SX_m$ is also monotone increasing in $m$ and converges
(downward) to $p$ as $m \to -\infty$.  Our aim is to close the argument by showing that such dynamical behavior can occur only
when $B=0$ for the quantity in Condition B, and thus can only happen when Condition B fails.
In other words, when Condition B holds and so $B\ne 0$, then
the inequality \eqref{CP} is strict at $\theta_0=0$
and so crystallographic pinning occurs in the direction $\theta_0 = 0$.

Our approach is to consider the difference equation 
\begin{equation}
\begin{array}{lcl}
Y_{m+1} &\bb = &\bb 2Y_m - Y_{m-1} - (S + S^{-1} - 2I)Y_m\\
\\
&\bb &\bb + f(p+Y_m,a_+(0)) - f(p,a_+(0)),
\end{array}
\label{4diff}
\end{equation}
which is satisfied by both $Y_m=X_m-p$ and also by $Y_m=SX_{-m}-p$ whenever $X_m$ satisfies \eqref{diff}.
We shall show that when Condition A holds, then the system \eqref{4diff} possesses a two-dimensional
center manifold at the origin.
Moreover, this center manifold contains two monotone orbits which approach $0$ from below and above,
and toward which the orbits $X_m - p$ and $SX_{-m} - p$, respectively, obtained in Proposition 2.4 converge exponentially fast
as $m\to\infty$.  A convexity argument then shows that
the existence of one of these orbits forces $B \le 0$ and the other forces $B \ge 0$, thus $B = 0$ must hold.

The first step in this process is the construction of a smooth center manifold.
The following theorem is a discrete-time version of Theorem 4.1 in \cite{Vanderbauwhede2}.
For details see \cite{James} in the case where $\X$ is a Hilbert space.
\vv

\noindent {\bf Theorem 4.1 (Center Manifold Theorem).} \bxx
Let $\X$ be a Banach space and let $A \in \mathcal{L}(\X)$ satisfy the following spectral gap condition:
There exists $\alpha > 0$ such that for each $\lambda \in \sigma(A)$ either
$|\lambda| = 1$ or $| \log |\lambda| |  > \alpha$.
Denote by $\X^c$ and $\X^h$ the associated center and hyperbolic subspaces of $\X$, corresponding
to spectra with $|\lambda|=1$ and $|\lambda|\ne 1$, respectively.
Assume that $\X^c$ is finite dimensional.  Let $g:\X\to\X$ be a $C^k$-smooth function
for some $1\le k<\infty$ and satisfy $g(0) = 0$ and $Dg(0) = 0$.  Then there exist
neighborhoods $\Omega^c\subseteq\X^c$ and $\Omega^h\subseteq\X^h$ of the origin in these
subspaces, and a
$C^k$ mapping $\psi:\Omega^c\to\Omega^h$ with $\psi(0) = 0$
and $D\psi(0) = 0$, such that the following properties hold:

\begin{itemize}
\item If $x^c_n\in\Omega^c$ for $0\le n\le N+1$ satisfies the reduced system
\begin{equation} x^c_{n+1} = A^c x_n^c + \pi^cg(x_n^c+\psi(x_n^c)) \label{reduced} \end{equation} for $0 \le n \le N$,
and if we let 
$x_n = x_n^c + \psi(x_n^c)$, then $x_n$ satisfies the full system
\begin{equation} x_{n+1} = Ax_n + g(x_n) \label{full} \end{equation}  for $0 \le n \le N$.
\item If $x_n$ satisfies the full system
\eqref{full} and also $x_n \in \Omega$ for all $n \in \Z$,
then $\pi^hx_n = \psi(\pi^c x_n)$
and $x_n^c=\pi^c x_n$ satisfies the reduced system \eqref{reduced}
for every $n\in\Z$.
\end{itemize}
In the above statements, we let $\pi^c$ and $\pi^h$ denote the canonical projections with respect to
the decomposition $\X=\X^c\oplus\X^h$,
and $A^c=\pi^cA$ and $A^h=\pi^hA$ are the corresponding operators on these subspaces. We also denote
$\Omega=\{x^c+x^h\;|\;x^c\in\Omega^c\hbox{ and }x^h\in\Omega^h\}$, which is a neighborhood of the origin in $\X$.
\exx
\vv

When Theorem 4.1 holds, we call the set
$$
W^c = \{x^c + \psi(x^c)\; | \; x^c \in \Omega^c\}\subseteq\Omega
$$
a {\bf local center manifold} for the full system \eqref{full}, and the
system \eqref{reduced} restricted to this set is called the {\bf reduced system}.
We only outline the proof of this result. For more details, see \cite{James,Vanderbauwhede1, Vanderbauwhede2}.
\vv

\noindent {\bf Sketch of Proof.}
The full evolution equation \eqref{full} may be written in terms of its center, stable, and unstable parts,
where we decompose $\X=\X^c\oplus\X^h=\X^c\oplus(\X^s\oplus\X^u)$ in this fashion. We generally
denote $x^\circ=\pi^\circ x$, with $A^\circ=\pi^\circ A$ and $g^\circ(x)=\pi^\circ g(x)$,
where $\pi^\circ$ denotes the spectral projection onto $\X^\circ$, with $\circ$ representing
$c$, $h$, $s$, or $u$. We thus have
$$
\ba{lcl}
x^c_{n+1} &\bb = &\bb A^c x^c_n + g^c(x_n),\\
\\ 
x^s_{n+1} &\bb = &\bb A^s x^s_n + g^s(x_n),\\
\\
x^u_{n+1} &\bb = &\bb A^u x^u_n + g^u(x_n),
\ea
$$
for the full equation \eqref{full}, and we have that
$x_n=x^c_n+x^s_n+s^u_n$.
The corresponding variation of constants formulae are
\begin{equation}
\ba{lcl}
x^c_n &\bb = &\bb \ds{(A^c)^nx^c_0 + \sum_{j=0}^{n-1} (A^c)^{n-1-j} g^c(x_j),}\\
\\
x^s_n &\bb = &\bb \ds{(A^s)^nx^s_0 + \sum_{j=0}^{n-1} (A^s)^{n-1-j} g^s(x_j),}\\
\\
x^u_n &\bb = &\bb \ds{(A^u)^nx^u_0 + \sum_{j=0}^{n-1} (A^u)^{n-1-j} g^u(x_j),}
\ea
\label{ev2}
\end{equation}
where if $n\le 0$ we interpret $\sum_{j=0}^{n-1}=-\sum_{j=n}^{-1}$, with $\sum_{j=0}^{-1}$ the empty sum.
We begin the proof by restricting attention to the case where the orbit $\{x_n\}_{n \in \Z}$
is bounded for all $n \in \Z$. In this case, upon multiplying the second (stable) equation in \eqref{ev2}
by $(A^s)^{-n}$ and letting $n\to-\infty$, and also multiplying the third (unstable) equation in \eqref{ev2}
by $(A^u)^{-n}$ and letting $n\to\infty$, we obtain
\be
\ba{lcl}
x^s_0 &\bb = &\bb \ds{\sum_{j=-\infty}^{-1} (A^s)^{-1-j} g^s(x_j),}\\
\\
x^u_0 &\bb = &\bb \ds{-\sum_{j=0}^\infty (A^u)^{-1-j} g^u(x_j),}
\ea
\label{ev3}
\ee
and then substituting \eqref{ev3} into the second and third equations
of \eqref{ev2} gives
\be
\ba{lcl}
x^s_n &\bb = &\bb \ds{\sum_{j=-\infty}^{n-1} (A^s)^{n-1-j} g^s(x_j),}\\
\\
x^u_n &\bb = &\bb \ds{-\sum_{j=n}^\infty (A^u)^{n-1-j} g^u(x_j).}
\ea
\label{ev4}
\ee
Combining \eqref{ev4} with the first equation in \eqref{ev2} thus gives
\be
\ba{lcl}
x_n &\bb = &\bb \ds{(A^c)^nx^c_0 + \sum_{j=0}^{n-1} (A^c)^{n-1-j} g^c(x_j)}\\
\\
&\bb &\bb \ds{+ \sum_{j=-\infty}^{n-1} (A^s)^{n-1-j} g^s(x_j) 
- \sum_{j=n}^\infty (A^u)^{n-1-j} g^u(x_j).}
\ea
\label{CME1}
\ee
We regard \eqref{CME1} as a fixed point equation for the trajectory
$\{x_n\}_{n \in \Z}$ with $x^c_0$ as a parameter.  For each $\zeta > 0$
define the Banach space $Y_\zeta$ of sequences $x=\{x_n\}_{n\in\Z}$ in $\X$ by 
$$
\Y_\zeta = \{x:\Z\to\X \; | \; \|x\|_\zeta<\infty\},\qquad
\|x\|_\zeta=\sup_{ n \in \Z} e^{-\zeta |n|} \|x_n\|_\X,
$$
with $\|\cdot\|_\zeta$ being the norm in $\Y_\zeta$. Also define
$G : \X^c \times \Y_\alpha \to \Y_\alpha$ to be the right hand side of equation \eqref{CME1},
with arguments $x^c_0\in\X^c$ and $x=\{x_n\}_{n\in\Z}\in \Y_\alpha$, with $\alpha$ as in
the statement of the theorem. Thus the fixed points of $G(x^c_0,\cdot)$ are
solutions of \eqref{CME1}.  It is not hard to show that $G(x^c_0,\cdot)$ is a contraction mapping,
so long as the lipschitz constant of $g$ is sufficiently small.
Assuming this,
denote the unique fixed point of $G(x^c_0,\cdot)$ by $\Psi(x^c_0)\in \Y_\alpha$,
and let $\psi : \X^c \to \X^h$ be given by $\psi(x^c_0) = \pi^h [\Psi(x^c_0)_0]$.
Observe that $x^c_n=\pi^c[ \Psi(x^c_0)_n]$ solves the reduced system \eqref{reduced}
and that
$x_n=\pi^c[ \Psi(x^c_0)_n] + \psi( \pi^c[ \Psi(x^c_0)_n ])$
solves the full system \eqref{full}. 

The difficult part is proving that $\psi$ is smooth.  Since we are using the uniform contraction mapping principle, we obtain that the fixed point $\Psi$ is only as smooth in $x^c_0$ as the uniform contraction mapping $G$ is.  The Nemytskii operator associated to $g$ is in general not smooth as a mapping from $\Y_\alpha$ to itself.  However, it is $C^k$-smooth as a mapping from $\Y_\zeta$ to $\Y_{\alpha}$ whenever $k\zeta<\alpha$.  Using this fact,
one may show with a clever application of the fiber contraction theorem to the map $G$
that $\Psi:\X^c\to \Y_\alpha$ is $C^k$-smooth, and thus $\psi:\X^c\to\X^h$ is $C^k$-smooth.

To complete the theorem, it remains only to remove the assumption that
the lipschitz constant of $g$ is small.  This is done by considering the modified system 
\begin{equation}
x_{n+1} = Ax_n + \tilde{g}(x_n)
\label{modified}
\end{equation}
where $\tilde{g}(x) = g(x)\chi(x)$.
Here $\chi(x) = \chi_0(\frac{\|x^c\|}{\eps}) \chi_0(\frac{\|x^h\|}{\eps})$
where $\chi_0 :[0,\infty)\to \R$ is a smooth cutoff function which is identically $1$ on $[0,1]$ and vanishes identically on $[2,\infty)$.  By choosing $\eps$ small enough, we can make
the lipschitz constant of $g$ as  small as we wish on the region
on which $\chi(x)$ is nonzero, and it follows from what we have sketched above that the system
\eqref{modified} has a global center manifold.
As for smoothness of this manifold, recall that in general Banach spaces need not have smooth norms.
However, as we have assumed that $\X^c$ is finite dimensional, we may assume that its norm is smooth.
Thus $\chi(x)$ depends smoothly on $x^c$, but in general is only lipschitz in $x^h$. It follows that
$\chi$ and thus $G(x^c_0,\cdot)$ is smooth in the region $\|x^h\|<\eps$ where the cutoff function in the hyperbolic
direction is constant.
One now uses exponential estimates on the map $G$ to show that
$\|\Psi(x^c_0)_n\|<\eps$ for all coordinates of the fixed point $\Psi(x^c_0)$ of $G(x^c_0,\cdot)$,
whenever $\|x^c_0\|<\eps$.  With this and the above remarks about smoothness, it follows that
$\psi(x^c_0)$ depends smoothly on $x^c_0$ for $\|x^c_0\|<\eps$. Taking $\Omega^c$ and $\Omega^h$ to
be the balls of radius $\eps$ in $\X^c$ and $\X^h$ centered at the origin, the result now follows.
\qed
\vv

In addition to the existence of a smooth local center manifold, our proof of Theorem 1.1 requires the following shadowing lemma which guarantees that each orbit which stays close to the center manifold approaches an orbit on the center manifold exponentially fast.  This is classical, and generally a consequence of the existence of a center-stable foliation.  We do not require the full apparatus here so we provide an independent proof for the specific estimate that we need.  The proof is similar to the construction of a stable manifold.
\vv

\noindent {\bf Lemma 4.2 (Shadowing Lemma).} \bxx
Consider the setting of Theorem 4.1, with $\alpha>0$
as in the statement of that result. Then there exists a neighborhood
$\Omega\subseteq\X$ of the origin such that any
forward solution $x=\{x_n\}_{n\ge 0}\subseteq\Omega$ to the full system \eqref{full}
which lies in $\Omega$ possesses
an exponentially close shadow on the center manifold.
More precisely, there exists a positive constant $K$ such that
for any such solution to equation \eqref{full}, there exists a sequence
$y^c = \{y^c_n\}_{n\ge 0} \subseteq \X^c$ which satisfies the reduced system \eqref{reduced}
and such that 
$$
\|x_n-y^c_n - \psi(y^c_n)\| \le K e^{-\alpha n} \|x_0-x^c_0 - \psi(x^c_0)\|,  \qquad n \ge 0.
$$
\exx
\vv

\noindent {\bf Proof.}
The proof breaks up naturally into two parts. First we show that $x_n$
approaches the center manifold exponentially fast. Then we show that $x_n$
differs from a particular orbit on the center manifold by an exponentially
decreasing amount. Throughout, we work in a sufficiently small neighborhood
$\Omega$ of the origin in $\X$.
In what follows we shall let $\tilde\Omega$ denote
the neighborhood which was denoted by $\Omega$ in the statement of Theorem 4.1,
and which contains the center manifold. We shall let $\Omega$ denote
a (possibly) smaller neighborhood $\Omega\subseteq\tilde\Omega$
which will be constructed below and
for which the statement of the present lemma is valid.
We shall otherwise keep the same notation as in the proof of Theorem 4.1.

We first introduce new coordinates $(x^c,w)$ replacing $(x^c,x^h)$
in a neighborhood of the origin, where $w\in\X^h$ replaces $x^h$ and
is defined by $w = x^h - \psi(x^c)$, and with the coordinate $x^c$ unchanged.
Thus $w=0$ if and only if $x$ lies on the center manifold, in this neighborhood,
and we have $x=x^c+\psi(x^c)+w$.
The evolution equation for the full system \eqref{full}
written in the new coordinates takes the form
\begin{equation}
\ba{lcl}
x^c_{n+1} &\bb = &\bb A^c x^c_n + g^c(x^c_n+\psi(x^c_n)+w_n), \\ \\
w_{n+1} &\bb = &\bb A^h w_n + A^h\psi(x^c_n)
+g^h(x^c_n +\psi(x^c_n)+w_n)\\ \\
&\bb &\bb - \psi(A^c x^c_n + g^c(x^c_n + \psi(x^c_n)+w_n)).
\ea
\label{SL2}
\end{equation}
It is a consequence of the invariance of the center manifold that $\psi$ satisfies the functional relation
$$
\psi(A^c u + g^c(u + \psi(u))) = A^h \psi(u) + g^h(u + \psi(u)),
$$
and it follows from this that we may rewrite the system \eqref{SL2} as
$$\ba{lcl}
x^c_{n+1} &\bb = &\bb A^c x^c_n + \tilde g^c(x^c_n,w_n),\\ \\
w_{n+1} &\bb = &\bb A^h w_n + \tilde g^h(x^c_n,w_n),
\ea
$$
where $\tilde g^c$ and $\tilde g^h$ are given by
$$
\ba{lcl}
\tilde g^c(u,w) &\bb = &\bb g^c(u+\psi(u)+w),\\
\\
\tilde g^h(u,w) &\bb = &\bb
g^h(u + \psi(u) + w) - g^h(u + \psi(u))\\ \\
&\bb &\bb +\psi(A^c u + g^c(u + \psi(u))) - \psi(A^c u + g^c(u+\psi(u)+w)).
\ea
$$
Observe that $\tilde g^h(u,0)=0$ identically, and also that
$D\tilde g^h(0,0)=0$ since $Dg(0)=0$ and $D\psi(0)=0$.


Now let $w^s=\pi^sw$ and $w^u=\pi^uw$ denote the projections of $w$ onto
the stable and unstable subspaces, and so $w=w^s+w^u$. Let
$\tilde g^s(u,w)=\pi^s\tilde g^h(u,w)$ and
$\tilde g^u(u,w)=\pi^u\tilde g^h(u,w)$,
and note that $\tilde g^s(u,0)=0$ and $\tilde g^u(u,0)=0$ identically,
and that $D\tilde g^s(0,0)=0$ and $D\tilde g^u(0,0)=0$.
Much as in the proof of Theorem 4.1,
we may write the variation of constants formulae
\be
\ba{lcl}
w^s_n &\bb = &\bb
\ds{(A^s)^n w^s_0 + \sum_{j = 0}^{n-1} (A^s)^{n-1-j}\tilde{g}^s(x^c_j,w_j),} \\ \\
w^u_n &\bb = &\bb
\ds{-\sum_{j = n}^\infty (A^u)^{n-1-j}\tilde g^u(x^c_j,w_j),}
\ea
\label{wsu}
\ee
although we take a finite sum here for the stable part in
contrast to the proof of Theorem 4.1. (We have also omitted
the formula for $x^c_n$ as it will not be needed here.)
Adding the two formulae in \eqref{wsu} gives
\be
w_n=(A^s)^nw^s_0+\sum_{j=0}^{n-1}(A^s)^{n-1-j}\tilde g^s(x^c_j,w_j)
-\sum_{j=n}^\infty (A^u)^{n-1-j}\tilde g^u(x^c_j,w_j),
\label{fp}
\ee
for every $n\ge 0$, and which is valid for any forward orbit $x=\{x_n\}_{n\ge 0}$
lying in $\tilde\Omega$.

We wish to show that
the coordinates $w_n$ decay exponentially to zero provided the sequence $x_n$
lies in a sufficiently small neighborhood $\Omega$
of the origin. Our approach is to
regard the coordinates $x^c_n$ and also $w^s_0$ as known, and to consider
equation \eqref{fp} as a fixed point problem for the bounded sequence $w=\{w_n\}_{n\ge 0}$.
More precisely, we shall show this sequence is in fact a fixed point
both in a space of bounded sequences, and also in
a space of exponentially decaying sequences.

To this end, for any $\zeta\ge 0$
define the Banach space $\ZZ_\zeta$ of one-sided sequences in~$\X^h$
$$
\ZZ_\zeta=\{z:\N_0\to\X^h\;|\;\|z\|_\zeta<\infty\},\qquad
\|z\|_\zeta=\sup_{n\ge 0}e^{\zeta n}\|z_n\|_\X,
$$
where $\N_0=\{0,1,2,\ldots\}$. (Contrast the definition of this space
of one-sided decaying sequences, with that of $\Y_\zeta$ which has
two-sided growing sequences. We trust that the same notation $\|\cdot\|_\zeta$
for the two norms will not lead to confusion.)
Next, with $\alpha$ as in the statement of the theorem, there exist
$\gamma>\alpha$ and $C>0$ such that
\be
\|(A^s)^n\|,\|(A^u)^{-n}\|\le Ce^{-\gamma n},\qquad n\ge 0.
\label{expbnd}
\ee
Now fix a neighborhood $\Omega\subseteq\tilde\Omega$ of the origin
and take small enough quantities $\eps$ and $\delta$
such that the following all hold. First,
\be
\|w^s\|\le\frac{\eps}{2C}\hbox{ whenever }x\in\Omega,
\label{wbnd}
\ee
where as usual $w^s=\pi^sw$ with $w=x^h-\psi(x^c)$. Next,
\be
\begin{array}{l}
\|D_z\tilde g^s(x^c,z)\|,\|D_z\tilde g^u(x^c,z)\|\le\delta\\
\\
\hbox{whenever }x\in\Omega\hbox{ and }z\in\X^h\hbox{ with }\|z\|\le\eps,
\end{array}
\label{lipg}
\ee
with $D_z$ denoting the derivative with respect to the second argument.
Here $z$ need not be the $w$-coordinate arising from $x$ as above,
but rather is an arbitrary point in the closed $\eps$-ball about the
origin in $\X^h$.
Finally, we require that the inequality \eqref{contrac2} below
with \eqref{contrac1} should hold
both for $\zeta=\alpha$ and for $\zeta=0$.
It is easily seen that all this can be accomplished by first fixing
$\delta$, and then $\eps$, and then taking
$\Omega$ sufficiently small.

Now let $\{x_n\}_{n\ge 0}\subseteq\Omega$ be as in the statement of the lemma,
with $w_n=x^h_n-\psi(x^c_n)$ as usual. Denote the closed $\eps$-ball
about the origin in $\ZZ_\zeta$ by
$$
B_\zeta(\eps)=\{z=\{z_n\}_{n\ge 0}\in \ZZ_\zeta\;|\;\|z\|_\zeta\le\eps\},
$$
and define a map $G:B_\zeta(\eps)\to \ZZ_\zeta$ by setting
$$
G(z)_n=(A^s)^nw^s_0+\sum_{j=0}^{n-1}(A^s)^{n-1-j}\tilde g^s(x^c_j,z_j)
-\sum_{j=n}^\infty (A^u)^{n-1-j}\tilde g^u(x^c_j,z_j)
$$
for $n\ge 0$. Of course it must be shown that $G(z)$ actually lies
in $\ZZ_\zeta$. We show this, and will also bound the lipschitz constant
of the map $G$, which will give conditions under which $G$ is a contraction mapping.

Fix $\zeta$ satisfying $0\le\zeta<\gamma$ and take any
$z,\overline z\in B_\zeta(\eps)$.
Then $\|z_n-\overline z_n\|\le e^{-\zeta n}\|z-\overline z\|_\zeta$, and
using \eqref{expbnd} and \eqref{lipg} we have that
$$
\ba{l}
\|G(z)_n-G(\overline z)_n\|\\ \\
\qquad \ds{\le C\delta\sum_{j=0}^{n-1}e^{-\gamma(n-1-j)}e^{-\zeta j}\|z-\overline z\|_\zeta
+C\delta\sum_{j=n}^\infty e^{-\gamma(j+1-n)}e^{-\zeta j}\|z-\overline z\|_\zeta}\\ \\
\qquad < \tilde C(\zeta)C\delta e^{-\zeta n}\|z-\overline z\|_\zeta
\ea
$$
after a short calculation, where
\be
\tilde C(\zeta)=\frac{e^\zeta}{1-e^{-(\gamma-\zeta)}}+\frac{e^{-\gamma}}{1-e^{-(\gamma+\zeta)}}.
\label{contrac1}
\ee
Thus we have that
$$
\|G(z)-G(\overline z)\|_\zeta\le\tilde C(\zeta)C\delta\|z-\overline z\|_\zeta.
$$
Next note that the sequence $(A^s)^nw^s_0$, which is $G(0)$, lies in the space $\ZZ_\zeta$
with the bound
\be
\|G(0)\|_\zeta\le C\|w^s_0\|\le\frac{\eps}{2}
\label{zerbnd}
\ee
on its norm, from \eqref{expbnd} and \eqref{wbnd}.
It follows immediately that if
\be
\tilde C(\zeta)C\delta\le\frac{1}{2}
\label{contrac2}
\ee
holds, then $G$ is a contraction mapping of $B_\zeta(\eps)$ into itself, and thus
has a fixed point. In fact, we have assumed this to be the case for the choices
$\zeta=0$ and $\zeta=\alpha$. Moreover, we have that $B_\alpha(\eps)\subseteq B_0(\eps)$
and so the fixed points in these two balls are identical, namely, the bounded sequence
$\{w_n\}_{n\ge 0}$ obtained from the original orbit $\{x_n\}_{n\ge 0}$ and
which satisfies \eqref{fp}.
Thus $\{w_n\}_{n\ge 0}$ in fact lies in $B_\alpha(\eps)$.
Additionally, as the contraction constant of $G$ in this ball is bounded by $\frac{1}{2}$, it
follows that $\|w\|_\alpha\le 2\|G(0)\|_\alpha$, which with \eqref{zerbnd} yields the bound
$$
\|w_n\|\le 2Ce^{-\alpha n}\|w^s_0\|.
$$

We now show the existence of a shadow orbit $\{y^c_n\}_{n\ge 0}$, as in
the statement of the lemma. As this part of the proof is very similar to the
part above, and in particular involves a contraction mapping in the space
$\ZZ_\alpha$, we only sketch the argument.
We write $u_n=y^c_n-x^c_n$. As $x^c_n$ satisfies the first equation in \eqref{SL2}
and we wish $y^c_n$ to satisfy the reduced equation \eqref{reduced}, then $u_n$ must satisfy
\be
u_{n+1}=A^cu_n+\hat g(x^c_n,u_n)+q_n
\label{ueq}
\ee
where
$$
\hat g(x^c,u)=\tilde g^c(x^c+u,0)-\tilde g^c(x^c,0),\qquad
q_n=\tilde g^c(x^c,0)-\tilde g^c(x^c,w_n).
$$
A sufficient condition for \eqref{ueq} to hold for $n\ge 0$,
with $u_n$ decaying exponentially, is that
\be
u_n=Q_n-\sum_{j=n}^\infty(A^c)^{n-1-j}\hat g(x^c_j,u_j)
\label{fpu}
\ee
where
$$
Q_n=-\sum_{j=n}^\infty(A^c)^{n-1-j}q_j.
$$
In particular, one has the (growing) estimate
$$
\|(A^s)^{-n}\|\le Ce^{\kappa n},\qquad n\ge 0,
$$
where $\kappa>0$ can be taken arbitrarily small, and where
$C=C(\kappa)$. One checks that $\{q_n\}_{n\ge 0}$ and thus
$\{Q_n\}_{n\ge 0}$ belong to the space $\ZZ_\alpha$, and also
that equation \eqref{fpu} possesses a fixed point near the
origin in this space.  This concludes the sketch of the proof of the existence of a shadow orbit and hence the proof of the lemma.
\qed
\vv

The existence of a center manifold for our particular system \eqref{diff} follows from the center manifold theorem upon checking that we have a spectral gap at the unit circle.  The following lemma establishes the existence of a spectral gap and in addition characterizes the center subspace for a class of difference equations which includes \eqref{diff}. 
\vv

\noindent {\bf Lemma 4.3.} \bxx
Let $L\in\mathcal{L}(\X)$ be an operator for which $0\in\sigma(L)$ is an isolated
point of the spectrum, and let $\pi^0\in\mathcal{L}(\X)$
denote the spectral projection onto the
spectral subspace of $\X$ corresponding to this point.
Let $\ZZ=\X\times\X$ and define $A\in\mathcal{L}(\ZZ)$ by
$$
A = \left(\ba{cc} 2I-L & -I \\ I & 0\strt \ea \right).
$$
Then $1\in\sigma(A)$ is an isolated point in the spectrum of $A$, and
the spectral projection $\Pi\in\mathcal{L}(\ZZ)$ corresponding to this point is given by
\be
\Pi=\left(\ba{cc} \pi^0 & 0 \\ 0 & \pi^0\strt\ea\right).
\label{bigpi}
\ee
Also, if there exists $\beta_1>0$ such that
\be
\sigma(L)\subseteq\{\lambda\in\C\;|\;\re\lambda<-\beta_1\}\cup\{0\}
\label{specleft}
\ee
then there exists $\beta_2>0$ such that
\be
\sigma(A)\subseteq\{\lambda\in\C\;|\;| |\lambda| - 1| > \beta_2\}\cup\{1\}.
\label{beta2}
\ee
Finally, if the operator $L$ is Fredholm of index zero and has a one-dimensional
generalized kernel spanned by $v\in\X\setminus\{0\}$, that is
$\ker(L)=\{av\;|\;a\in\R\}$ with $v\not\in\ran(L)$,
then
$$
\ran(\pi^0)=\{av\;|\;a\in\R\},\qquad
\ran(\Pi)=\{\col(av,bv) \; | \; (a,b) \in \R^2\},
$$
are the above-mentioned spectral subspaces of $L$ and $A$, and they
have dimension one and two, respectively.
\exx
\vv

\noindent {\bf Proof.}
We first show that
\be
\sigma(A)\subseteq\{\lambda\in\C\setminus\{0\}\;|\;2-\lambda-\lambda^{-1}\in\sigma(L)\}.
\label{speca}
\ee
Suppose that $\lambda\ne 0$ is such that $2-\lambda-\lambda^{-1}\not\in\sigma(L)$,
and denote $\Gamma(\lambda)=((2-\lambda-\lambda^{-1})I-L)^{-1}$. We claim that
\be
(\lambda I-A)^{-1}=\left(\ba{cc}
-\lambda I & I\\
-I & (2-\lambda)I-L\strt
\ea\right)
\left(\ba{cc}
\lambda^{-1}\Gamma(\lambda) & 0\\
0 & \lambda^{-1}\Gamma(\lambda)\strt
\ea\right),
\label{resolv}
\ee
and thus $\lambda\not\in\sigma(A)$.
Indeed, this is straightforwardly
proved by multiplying the above matrix product by $\lambda I-A$,
to obtain the identity matrix after a brief calculation. We omit the details.
Thus \eqref{speca} follows from this, and from the fact that $0\not\in\sigma(A)$,
which holds because
$$
A^{-1}=\left(\ba{cc}
0 & I\\
-I & 2I-L\strt
\ea\right).
$$
It now follows from \eqref{speca}, and from the easily-checked fact that
$2-\lambda-\lambda^{-1}\ne 0$ whenever $\lambda\ne 1$, that if $\lambda=1$
belongs to the spectrum of $A$ then it is an isolated point of the spectrum.

It is also easy to establish \eqref{beta2} assuming that \eqref{specleft}
holds. Indeed, with \eqref{specleft} holding
assume that $\lambda\in\sigma(A)\setminus\{1\}$. Then
$2-\lambda-\lambda^{-1}\in\sigma(L)\setminus\{0\}$ and so
$2-\re(\lambda+\lambda^{-1})<-\beta_1$. Writing $\lambda=|\lambda|e^{i\theta}$, we
see that
$$
|\lambda|+|\lambda|^{-1}
\ge (|\lambda|+|\lambda|^{-1})\cos\theta=\re(\lambda+\lambda^{-1})>2+\beta_1.
$$
Denoting the two positive roots of $r+r^{-1}=2+\beta_1$ by $r_\pm$ with $r_-<1<r_+$,
we have that either $|\lambda|<r_-$ or $|\lambda|>r_+$, and thus \eqref{beta2} holds
with $\beta_2=\min\{1-r_1,r_+-1\}$.

We next calculate the spectral projection $\Pi$ for $A$ corresponding to
the point $\lambda=1$. (At this point we do not yet know that $1\in\sigma(A)$,
however, this fact will follow when we show that $\Pi\ne 0$.) For sufficiently
small $r$ we have that
$$
\Pi=\frac{1}{2\pi i}\int_{|\lambda-1|=r}(\lambda I-A)^{-1}\:d\lambda.
$$
We calculate the above integral for each of the four block entries of
the matrix \eqref{resolv}, but for simplicity we shall only provide the
details for one of the entries, as the approach for the others is similar.
We take the lower right-hand entry
\be
\ba{lcl}
\Pi_{2,2} &\bb = &\bb
\ds{\frac{1}{2\pi i}\int_{|\lambda-1|=r}((2-\lambda)I-L)\lambda^{-1}\Gamma(\lambda)\:d\lambda}\\
\\
&\bb = &\bb
\ds{\frac{1}{2\pi i}\int_{|\nu|=r}(L+(\nu-1)I)(\nu^2I+(\nu+1)L)^{-1}\:d\nu,}
\ea
\label{contour}
\ee
with $\nu=\lambda-1$,
as one sees after a short calculation using the above formula for $\Gamma(\lambda)$.
We apply the projections $\pi^0$ and $\pi^1$, where we denote $\pi^1=I-\pi^0$, to the second integral in \eqref{contour}
to get
\be
\pi^j\Pi_{2,2}=\frac{1}{2\pi i}\int_{|\nu|=r}(L^j+(\nu-1)I)(\nu^2I+(\nu+1)L^j)^{-1}\:d\nu
\label{projint}
\ee
for $j=1,2$, where $L^j=\pi^j L$ is regarded as an operator on the spectral
subspace $\pi^j\X\subseteq\X$, and we separately calculate these contour integrals.
Now $0\not\in\sigma(L^1)$, and thus the integrand
$(L^1+(\nu-1)I)(\nu^2I+(\nu+1)L^1)^{-1}$ in \eqref{projint}
is holomorphic in $\nu$
in a neighborhood of $\nu=0$. Thus the integral vanishes by Cauchy's theorem,
and so $\pi^1\Pi_{2,2}=0$. On the other hand,
for the case $j=0$ we have $\sigma(L^0)=\{0\}$,
and so the integrand is holomorphic for all complex $\nu$ except $\nu=0$. In this case we may
increase the radius $r$ of the contour arbitrarily, again by Cauchy's theorem.
Upon doing so, and then scaling $\nu$, we
write the integral in the equivalent form
$$
\pi^0\Pi_{2,2}=\frac{1}{2\pi i}\int_{|\nu|=1}\nu^{-1}((r\nu)^{-1}L^0+(1-(r\nu)^{-1})I)
(I+((r\nu)^{-1}+(r\nu)^{-2})L^0)^{-1}\:d\nu.
$$
In the limit $r\to\infty$ we obtain
$$
\pi^0\Pi_{2,2}=\frac{1}{2\pi i}\int_{|\nu|=1}\nu^{-1}I\:d\nu=I,
$$
which is the identity operator on the space $\pi^0\X$. It thus follows
that $\Pi_{2,2}=\pi^0$. The proofs for the remaining three cases are similar, in which one shows
that $\Pi_{1,1}=\pi^0$ and $\Pi_{1,2}=\Pi_{2,1}=0$.
This establishes \eqref{bigpi}.

The final statement of the lemma follows immediately from the
fact that with $L$ having the Fredholm properties described,
we have $\X=\ker(L)\oplus\ran(L)$, and so
the range of the spectral projection $\pi^0$ is the one-dimensional
subspace $\ker(L)\subseteq\X$ spanned by $v$.
\qed
\vv

In the following proposition we construct a center manifold for the difference equation
\eqref{4diff} and compute the reduced equations on it. To begin, let us rewrite \eqref{4diff}
as a system in the product space $\ZZ=\X\times\X$ where $\X=\ell^\infty(\Z)$, namely
\be
\ba{lcl}
Y_{m+1} &\bb = &\bb (2I-L)Y_m-W_m+g(Y_m),\\
\\
W_{m+1} &\bb = &\bb Y_m,
\ea
\label{wysys}
\ee
where $L=L_p(a_+(0))$ is the operator in \eqref{lop} and where
$g:\X\to\X$ is given by
\be
g(Y)=f(p+Y,a_+(0))-f(p,a_+(0))-f'(p,a_+(0))Y.
\label{gform}
\ee
Thus we have the system $Z_{m+1}=AZ_m+G(Z_m)$
with $Z\in\ZZ$, and $A\in\mathcal{L}(\ZZ)$ and $G:\ZZ\to\ZZ$ given by
$$
Z=\left(\ba{c}Y\\ W\strt\ea\right),\qquad
A=\left(\ba{cc} 2I-L & -I \\ I & 0 \strt\ea\right),\qquad
G(Z)=\left(\ba{c}g(Y)\\ 0\strt\ea\right).
$$
The following center manifold reduction holds for this system.
\vv


\noindent {\bf Proposition 4.4.} \bxx
Let Condition A hold and take $v$ as in Proposition 1.4.
Then the system \eqref{wysys}, \eqref{gform} in the space $\ZZ$
satisfies the conditions of Theorem 4.1
(the Center Manifold Theorem)
with the associated center and hyperbolic subspaces $\ZZ^c$ and $\ZZ^h$
of $\ZZ$ given by
\be
\ba{lclcl}
\ZZ^c &\bb = &\bb \X^c\times\X^c, &\bb \qquad &\bb \X^c=\{av\;|\;a\in\R\},\\
\\
\ZZ^h &\bb = &\bb \X^h\times\X^h, &\bb \qquad &\bb
\X^h=\{x\in\X\;|\;\langle v,x\rangle=0\},
\ea
\label{ych}
\ee
and with a center manifold $W^c\subseteq\ZZ$ given by
\be
W^c=\{\col(\eta v+\psi_1(\eta,\omega),\omega v+\psi_2(\eta,\omega))\;|\;
(\eta,\omega)\in\Omega^c\subseteq\R^2\}
\label{wycent}
\ee
for smooth functions $\psi_i:\Omega^c\to\X^h$ satisfying $\psi_i(0)=0$
and $D\psi_i(0)=0$, where $\Omega^c\subseteq\R^2$ is a neighborhood of the origin.
The reduced system on $W^c$, in the coordinates $(\eta,\omega)\in\Omega^c$,
has the form
\be
\ba{lcl}
\eta_{m+1} &\bb = &\bb 2\eta_m-\omega_m+\gamma(\eta_m,\omega_m),\\
\\
\omega_{m+1} &\bb = &\bb \eta_m,
\ea
\label{red2}
\ee
where $\gamma:\Omega^c\to\R$ satisfies
\be
\gamma(\eta,\omega)=\langle v,g(\eta v+\psi_1(\eta,\omega))\rangle
=B\eta^2+o(|\eta|^{2}+|\omega|^{2}),
\label{gamma2}
\ee
where $B\in\R$ is the quantity \eqref{bsum} in Condition B\@.
\exx
\vv

%
%
%
%
%

\noindent {\bf Proof.}
Proposition 1.4 implies that
the operator $L=L_p(a_+(0))$ satisfies all the conditions of Lemma 4.3,
including \eqref{specleft}, along with the Fredholm condition and statements about
its kernel and range, and with the spectral projection $\pi^0\in\LL(\X)$
given by $\pi^0 x=\langle v,x\rangle v$.
It thus follows by Lemma 4.3 and because $G(0)=0$ and $DG(0)=0$,
that the conditions of Theorem 4.1 (the Center Manifold Theorem) hold for the system \eqref{wysys}
with $\ZZ^c$ and $\ZZ^c$ as in \eqref{ych}. Thus there exists a center manifold \eqref{wycent}
with functions $\psi_i$ as stated.

All that remains is to verify the form of the reduced system. We take
$$
Z_m=\left(\ba{c}Y_m\\
W_m\strt\ea\right)
=\left(\ba{c}\eta_m v+\psi_1(\eta_m,\omega_m)\\
\omega_m v+\psi_2(\eta_m,\omega_m)\strt\ea\right)
$$
in $W^c$, and similarly with $Z_{m+1}$, in the system \eqref{wysys}. We then
apply the projection $\Pi$ given by \eqref{bigpi}, that is, we take the
inner product of each equation with $v$. In doing so we note that
$\langle v,\psi_i(\eta,\omega)\rangle=0$ identically since the range of $\psi_i$
lies in $\X^h$, and also that $\langle v,x\rangle=0$
for every $x\in\ran(L)$, by Proposition 1.4. This directly gives \eqref{red2}
with the first equality in \eqref{gamma2}. The second equality
in \eqref{gamma2} follows from the fact that $g(0)=0$ and $Dg(0)=0$, and
also $\psi_1(0,0)=0$ and $D\psi_1(0,0)=0$, which implies that
$$
\langle v,g(\eta v+\psi_1(\eta,\omega))\rangle=
\frac{1}{2}\langle v,D^2g(0)(v,v)\rangle \eta^2+o(|\eta|^{2}+|\omega|^{2}),
$$
where $D^2g(0)(\cdot,\cdot)$ denotes the usual bilinear form of
the second derivative. Continuing, we have
$$
\langle v,D^2g(0)(v,v)\rangle=\langle v,f''(p,a_+(0))(v,v)\rangle
=\sum_{n=-\infty}^\infty f''(p_n,a_+(0))v_n^3=2B,
$$
to give the result as claimed. This completes the proof.
\qed

\section{The Proof of Theorem 1.1}
Theorem 1.1 now follows from a convexity argument applied to the reduced equations.  The only remaining preparatory results guarantee that the shadow orbits corresponding to $X_m - p$ and $SX_{-m} - p$ are not identically zero.
 
The following lemma, which is related to an exercise in Coddington and Levinson's classic text, will allow us to prove that $X_m$ decays only polynomially fast, hence its shadow orbit on the center manifold whose existence is guaranteed by Lemma 4.2 is not the zero orbit.
\vv

\noindent {\bf Lemma 5.1.} \bxx
Let $d_n$ for $n\ge 1$ be a sequence of real numbers which satisfies
\begin{equation}
\sum_{n=1}^\infty n^3|d_n| < \infty,
\label{BoundEqn1}
\end{equation}
and consider the difference equation
\begin{equation}
x_{n+1} = (2+d_n)x_n -  x_{n-1}.
\label{E43}
\end{equation}
Then given any initial condition $(x_0,x_1)\in\R^2$, there exist quantities
$(P,Q) \in \R^2$ such that
\be
x_n=Pn+Q+o(1)
\label{linlim}
\ee
as $n\to\infty$. Conversely, given any $(P,Q)\in\R^2$, there exists a unique
initial condition $(x_0,x_1)\in \R^2$ such that \eqref{linlim} holds.
In particular, if ${\ds{\lim_{n\to\infty}}}x_n=0$ then $x_n=0$ for every~$n$.
\exx
\vv

\noindent {\bf Proof.}
Let $\{x_n\}_{n\ge 0}$ satisfy \eqref{E43} and
denote $u_n = \col(x_{n+1},x_n)$ for $n\ge 0$. Let $U_n$ be the transition matrix 
$$
U_n=U+d_{n+1}\Gamma,\qquad
U=\left(\ba{cc} 2 & -1 \\
1 & 0 \strt\ea\right),\qquad
\Gamma=\left(\ba{cc} 1 & 0 \\
0 & 0 \strt\ea\right),
$$
and so $u_{n+1}= U_n u_n$. Note here that
\be
U^n=\left(\ba{cc} n+1 & -n \\
n & -n+1 \strt\ea\right)
\label{upower}
\ee
for every integer $n$, as is easily proved by induction.
Now let $v_n=U^{-n}u_n$ and observe that $v_n$ evolves according to the equation
$$
v_{n+1}=(I+R_n)v_n,\qquad
R_n=d_{n+1}U^{-(n+1)}\Gamma U^n
=d_{n+1}\left( \ba{cc}
-n(n+1)  & n^2  \\
-(n+1)^2  & n(n+1) \strt\ea \right),$$
by a simple calculation.
Next, if $n\ge m\ge 0$ let
\be
T_{n,m}=(I+R_{n-1})(I+R_{n-2})\cdots(I+R_m)
\label{tproduct}
\ee
with $T_{n,n}=I$, and so $v_n=T_{n,m}v_m$. Observe that
\be
\ba{l}
\|T_{n,m}\| \le \ds{\prod_{j=m}^{n-1}(1+\|R_j\|)\le\exp\bigg(\sum_{j=m}^{n-1}\|R_j\|\bigg),}\\
\\
\|T_{n,m}-I\| \le \ds{\bigg(\prod_{j=m}^{n-1}(1+\|R_j\|)\bigg)-1
\le\exp\bigg(\sum_{j=m}^{n-1}\|R_j\|\bigg)-1.}
\ea
\label{estimates}
\ee
The first inequality in the second line of \eqref{estimates} may require
a brief explanation. This inequality is obtained by first expanding the matrix product
in \eqref{tproduct}, then subtracting the term $I$, thereby obtaining a polynomial in the matrices $R_j$.
Next one takes the norm of this polynomial, and passes the norm across all the terms,
obtaining the same polynomial but now in the scalar quantities $\|R_j\|$. The polynomial
so obtained is the second term in the second line of \eqref{estimates}, as desired.
 
We have for some $C_1>0$ that $\|R_n\|\le C_1(1+n^2)|d_{n+1}|$ for every $n\ge 0$,
and so from \eqref{BoundEqn1} and \eqref{estimates}
there exists $C_2>0$ such that $\|T_{n,m}\|\le C_2$ for every $n\ge m\ge 0$.
Further,
$$
\ba{lcl}
\|T_{n,0}-T_{m,0}\| &\bb = &\bb \ds{\|(T_{n,m}-I)T_{m,0}\| \le
C_2\bigg(\exp\bigg(\sum_{j=m}^{n-1}\|R_j\|\bigg)-1\bigg)}\\
\\
& \le &
\ds{C_2\bigg(\exp\bigg(\sum_{j=m}^\infty C_1(1+j^2)|d_{j+1}|\bigg)-1\bigg)=r_m}
\ea
$$
with the above equation serving as the definition of the
quantity $r_m$. As ${\ds{\lim_{m\to\infty}}}r_m=0$,
it follows that $T_{n,0}$ is a Cauchy sequence of matrices, so
the limit $T_{\infty,0}={\ds{\lim_{n\to\infty}}}T_{n,0}$
exists. The vectors $v_n$ as well possess a limit
$$
v_\infty=\lim_{n\to\infty}v_n=T_{\infty,0}v_0=T_{\infty,0}u_0.
$$
Let us also note the estimate
$$
mr_m\le mC_3\sum_{j=m}^\infty j^2|d_{j+1}|\le C_3\sum_{j=m}^\infty j^3|d_j|
$$
for some $C_3>0$, which follows from the mean-value theorem, and which, with \eqref{BoundEqn1},
implies that ${\ds{\lim_{m\to\infty}}}mr_m=0$.
We thus have that
$$
u_n=U^n v_n=U^n v_\infty+U^n(v_n-v_\infty),
$$
and as $\|U^n\|\le C_4(n+1)$ for some $C_4>0$, by \eqref{upower}, we have that
$$
\begin{array}{lcl}
\|U^n(v_n-v_\infty)\| &\bb \le &\bb C_4(n+1)\|T_{n,0}-T_{\infty,0}\|\|u_0\|\\
\\
&\bb \le &\bb
C_4(n+1)r_n\|u_0\|\to 0,\qquad\hbox{as }n\to\infty.
\end{array}
$$
Thus $u_n=U^nv_\infty+o(1)$, and denoting $v_\infty=\col(A,B)$ it follows
directly from \eqref{upower} that the second coordinate of $u_n$,
which is $x_n$, has the form
$$
x_n=An-B(n-1)+o(1).
$$
This proves \eqref{linlim} with $P=A-B$ and $Q=B$.

The converse is proved more or less by following the above steps
in reverse. Namely, given $(P,Q)\in\R^2$, then let
$v_\infty=\col(A,B)$ where $A=P+Q$ and $B=Q$, and let
$u_0=T_{\infty,0}^{-1}v_\infty$. We note that the matrix $T_{\infty,0}$
is invertible as $\det(I+R_j)=1$ for every $j$, and thus $\det T_{n,m}=1$
for every $n\ge m\ge 0$, and thus $\det T_{\infty,0}=1$. The required
initial condition is thus $u_0=\col(x_1,x_0)$, and we see that it is unique.

The final sentence in the statement of the lemma follows in particular
from the uniqueness of $(x_0,x_1)$ for a given $(P,Q)$.
\qed
\vv

\noindent {\bf Proposition 5.2.} \bxx
Let $x_n$ for $n\ge 0$ be a sequence of numbers satisfying
\be
x_{n+1}-2x_n+x_{n-1}=Mx_n^2+\rho(x_{n-1},x_n,x_{n+1}),\qquad n\ge 1,
\label{551}
\ee
and
\be
x_n\ne 0\hbox{ for infinitely many }n,\qquad
\lim_{n\to\infty}x_n=0,
\label{552}
\ee
where $M\ne 0$ and where the function $\rho:\R^3\to\R$ is $C^{2}$ and satisfies
\be
\rho(\beta_-,\beta_0,\beta_+)=o(|\beta_-|^{2}+|\beta_0|^{2}+|\beta_+|^{2})
\label{55r}
\ee
at the origin. Then
$$
3Mx_n>Mx_{n-1}>Mx_n>0\hbox{ for all large }n.
$$
\exx
\vv

\noindent {\bf Remark.} We expect that in fact ${\ds{\lim_{n\to\infty}}}\frac{x_n}{x_{n-1}}=1$ should
hold in the above proposition, although do not need this fact.
\vv

\noindent {\bf Proof.}
Without loss, we shall assume that $M>0$, as
the case $M<0$ follows by considering the sequence $-x_n$ in place of $x_n$.

We first show that
\be
(x_{n-1},x_n,x_{n+1})\ne (0,0,0)\hbox{ for all large }n.
\label{556}
\ee
If \eqref{556} is false, then there are infinitely many $n$ for which
$x_{n-1}=x_n=x_{n+1}=0$ but
$x_{n+2}\ne 0$, in light of the first statement in \eqref{552}.
For such $n$, equation \eqref{551} with $n+1$ replacing $n$ takes the form
$x_{n+2}=\rho(0,0,x_{n+2})$, and so
\be
\frac{\rho(0,0,x_{n+2})}{x_{n+2}}=1.
\label{55f}
\ee
However, equation \eqref{55f} cannot hold for infinitely many $n$
in light of \eqref{55r} and the second statement in \eqref{552}.
With this contradiction \eqref{556} is proved.

We next show that for all sufficiently large $n$, the two inequalities
\be
x_{n+1}\le x_n,\qquad x_{n-1}\le x_n,
\label{558}
\ee
cannot simultaneously hold. Define a function $H_1:\R^3\to\R$ by
$$
H_1(\beta_+,\beta_0,\beta_-)=\beta_+^2+\beta_-^2+M\beta_0^2+\rho(\beta_0-\beta_-^2,\beta_0,\beta_0-\beta_+^2),
$$
and note that $H_1$ has a strict local minimum $H_1(0,0,0)=0$ at the origin,
by \eqref{55r}.
Now if $n$ is such that both inequalities in \eqref{558} hold,
let $\beta_\pm=(x_n-x_{n\pm 1})^{1/2}$ and observe that
equation \eqref{551} becomes
$H_1(\beta_+,x_n,\beta_-)=0$ for this $n$.
If $n$ is large enough this forces $\beta_\pm=x_n=0$ due
to the strict local minimum of $H_1$ and the limit in \eqref{552},
and thus $x_{n-1}=x_n=x_{n+1}=0$.
But this cannot happen for infinitely many $n$
by \eqref{556}. Thus the inequalities \eqref{558} cannot simultaneously hold
for arbitrarily large $n$.

It follows from the failure of \eqref{558} for all large $n$ that
the sequence $x_n$ is eventually monotone. In fact, either
\be
x_n>x_{n+1}>0\hbox{ for all large }n,
\label{dn}
\ee
or else
\be
x_n<x_{n+1}<0\hbox{ for all large }n.
\label{up}
\ee
We wish to prove \eqref{dn}, so let us assume to the contrary that \eqref{up} holds.
We claim that
\be
x_{n+1}-2x_n+x_{n-1}<0\hbox{ for infinitely many }n.
\label{55x}
\ee
If \eqref{55x} is false, then we have
that $x_{n+1}-x_n\ge x_n-x_{n-1}>0$ for all large $n$.
Thus there exists $\delta>0$ such that $x_{n+1}-x_n\ge\delta$
for all large $n$. However, this forces $x_n\to\infty$ as $n\to\infty$,
which is false. This establishes \eqref{55x}.

Still assuming \eqref{up}, for every $n$ for which the inequality in \eqref{55x}
holds let $\beta_\pm$ and $\beta_0$ be such that
$$
\beta_+=-x_{n+1},\qquad \beta_0=x_{n+1}-x_n,\qquad
\beta_-^2=-x_{n+1}+2x_n-x_{n-1}.
$$
We have that $\beta_+>0$ and $\beta_0>0$, and from equation \eqref{551} that
\be
\ba{lcl}
0 &\bb = &\bb
\beta_-^2+M(\beta_0+\beta_+)^2+\rho(-\beta_+-2\beta_0-\beta_-^2,-\beta_+-\beta_0,-\beta_+)\\
\\
&\bb > &\bb
\beta_-^2+M\beta_0^2+M\beta_+^2+\rho(-\beta_+-2\beta_0-\beta_-^2,-\beta_+-\beta_0,-\beta_+)\\
\\
&\bb = &\bb H_2(\beta_+,\beta_0,\beta_-),
\ea
\label{55h}
\ee
with the above formula serving as the definition of the function $H_2:\R^3\to\R$.
The function $H_2$ has a strict local minimum $H_2(0,0,0)=0$ at the origin by \eqref{55r},
and with \eqref{55h} this forces $\beta_\pm=\beta_0=0$ if $n$ is sufficiently large.
But this contradicts $\beta_+>0$ as noted above, and so
\eqref{dn} is established, as desired.

To complete the proof we must show that $3x_n\ge x_{n-1}$ for all large $n$. Upon
dividing equation \eqref{551} by $x_{n-1}$, we have that
$$
-\frac{2x_n}{x_{n-1}}+1<\frac{x_{n+1}-2x_n}{x_{n-1}}+1
=\frac{Mx_n^2+\rho(x_{n-1},x_n,x_{n+1})}{x_{n-1}}\to 0
$$
as $n\to\infty$, in light of the ordering $x_{n-1}>x_n>x_{n+1}>0$ established above.
Thus $1-\frac{2x_n}{x_{n-1}}<\frac{1}{3}$ for all large $n$, to give the result.
\qed
\vv

It is a consequence of Lemma 5.1 that the orbits $X_m$ guaranteed by Proposition 2.4 decay to $p$
as $m\to\infty$, or to $S^{-1}p$ as $m\to-\infty$, at a subexponential rate,
as the following result shows.
\vv

\noindent {\bf Proposition 5.3.} \bxx
Assume that the inequality \eqref{CP} is an equality at $\theta_0=0$.
Also assume that Condition A holds, with $p$ as stated there, and
let $X_m$ be the solution to \eqref{diff} guaranteed by
Proposition 2.4. Let $Y_m = X_m-p$.
Then
\begin{equation}
\sum_{m=0}^\infty m^3\|Y_m\| = \infty.
\label{ysum}
\end{equation}
The same conclusion holds if instead we let $Y_m=SX_{-m}-p$.
\exx
\vv

\noindent {\bf Proof.}
For definiteness we take $Y_m=X_m-p$. Also, we shall
denote the coordinates of $Y_m\in\ell^\infty(\Z)$ by $Y_m=\{y_{n,m}\}_{n\in\Z}$.
Let $v\in\ell^\infty(\Z)\setminus\{0\}$ be as in
Proposition 1.4, so of course $v\in\ell^1(\Z)$, and let
$y_m = \langle v,Y_m\rangle$. From the strict positivity $v_n>0$
of Proposition 1.4 and the strict ordering \eqref{2het} and the limits \eqref{xlims}
in Proposition 2.4, it follows that
\be
y_m<0,\qquad\lim_{m\to\infty}y_m=0.
\label{nonzero}
\ee
Also, $Y_m$ satisfies \eqref{4diff}, so it follows
using \eqref{lindiff} that $y_m$ satisfies
$$
y_{m+1} = (2+d_m)y_m - y_{m-1}
$$
where
$$
d_m = \frac{\langle v,f(p + Y_m, a_+(0)) - f(p,a_+(0)) - f'(p,a_+(0)) Y_m\rangle}{y_m}.
$$
Here we have used the fact that
$$
\begin{array}{lcl}
\langle v,(S + S^{-1} - 2I)Y_m\rangle
&\bb = &\bb \langle (S + S^{-1} - 2I)v, Y_m \rangle\\
\\
&\bb = &\bb \langle f'(p,a_+(0))v, Y_m\rangle
= \langle v,f'(p,a_+(0))Y_m\rangle.
\end{array}
$$
%
%
We note that the conclusion of Lemma 5.1 fails for the sequence $y_m$, in light of \eqref{nonzero},
and so necessarily 
\begin{equation}
\sum_{m = 0}^\infty m^3 |d_m| = \infty.
\label{asum}
\end{equation}
By the mean value theorem there are quantities $\eps_{m,n}\in[y_{m,n},0]$ such that
$$
\begin{array}{r}
f(p_n+y_{n,m},a_+(0)) - f(p_n,a_+(0)) - f'(p_n,a_+(0))y_{n,m}\\
\\
=f''(p_n+\eps_{n,m},a_+(0))y_{n,m}^2,
\end{array}
$$
and we note that
$p_n+\eps_{n,m}\in[p_n+y_{n,m},p_n]=[x_{n,m},p_n]\subseteq[-1,1]$.
Thus
\begin{equation}
\begin{array}{lcl}
|d_m| &\bb = &\bb
\ds{\frac{1}{|y_m|}
\bigg|\sum_{n=-\infty}^\infty v_nf''(p_n + \eps_{n,m},a_+(0)) y_{n,m}^2\bigg|}\\
\\
&\bb \le &\bb
\ds{\frac{K\|Y_m\|}{|y_m|}\sum_{n=-\infty}^\infty v_n |y_{n,m}|=K\|Y_m\|,}
\end{array}
\label{fineq}
\end{equation}
where $K>0$ is an upper bound for $|f''(u,a_+(0))|$ in the interval $[-1,1]$. Note in particular, in the final equality
in \eqref{fineq}, that we have used the fact that $v_n>0$ and $y_{n,m}<0$.
With this, the desired conclusion \eqref{ysum} follows from \eqref{asum} and \eqref{fineq}.
\qed
\vv


With these rates of convergence established, we may now compare the orbit $Y_m$ to its shadow on the center manifold.  
\vv


\noindent {\bf Proof of Theorem 1.1.}
Suppose that Condition A holds but that the inequality \eqref{CP} is in fact an equality, and so
crystallographic pinning does not occur in the direction $\theta_0 = 0$.
We will show that condition $B$ fails, namely, that $B = 0$ for the quantity
$B$ in formula \eqref{bsum}. Therefore assume to the contrary that $B\ne 0$.

With the above assumptions, let $X_m$ denote the solution of \eqref{diff} guaranteed by Proposition 2.4.
We first let $Y_m = X_m - p$, and so $Y_m$ satisfies \eqref{4diff}, equivalently,
$Z_m=\col(Y_m,Y_{m-1})$ satisfies the system \eqref{wysys}, \eqref{gform}. Note also that $Y_m$ approaches $0$ monotonically from below
as $m\to\infty$.
It is a consequence of Lemma 4.2 (the Shadowing Lemma) and Proposition 4.4
that there exists a sequence $\{\eta_m\}_{m\ge m_0}$ for some sufficiently large $m_0$, such that
\begin{equation}
\eta_{m+1} -2\eta_m+ \eta_{m-1} = \gamma(\eta_m,\eta_{m-1})
\label{CMdiff}
\end{equation}
for $m\ge m_0+1$, and such that
\be
\|V_m\| \le Ke^{-\alpha m},\qquad
V_m=Y_m-\eta_m v - \psi_1(\eta_m,\eta_{m-1}),
\label{S4E10}
\ee
for some positive constants $K$ and $\alpha$, where the above equality serves as the definition of $V_m\in\ell^\infty(\Z)$.
Thus the sequence $\eta_mv+\psi_1(\eta_m,\eta_{m-1})$ for large $m$ is
the shadow orbit on the center manifold to the given orbit $Y_m$,
with the function $\psi_1(\eta,\omega)$ as in Proposition 4.4, and where equation \eqref{CMdiff} is simply \eqref{red2} rewritten.

Now Proposition 5.3 implies that we have the divergent sum \eqref{ysum}. In light of the estimate \eqref{S4E10},
the analogous sum for $V_m$ converges, so it follows that
$$
\delta_m=\frac{\|V_m\|}{\|Y_m\|}\to 0\hbox{ for some subsequence }m=m_j\to\infty,
$$
where the above equality serves to define $\delta_m$.
Note further that there exist arbitrarily large $m$ such that $\eta_m\ne 0$.
Thus Proposition 5.2 applies to the sequence $\eta_m$, using the form \eqref{gamma2}
of the function $\gamma(\eta,\omega)$ and the fact that $B\ne 0$ is assumed, and we thereby conclude that
\be
3B\eta_m>B\eta_{m-1}>B\eta_m>0\hbox{ for all large }m.
\label{S54}
\ee
Dividing the equation in \eqref{S4E10} by $\eta_m$ and using the inequalities immediately above gives
\be
\lim_{m\to\infty}\frac{Y_m-V_m}{\eta_m}=v,
\label{S51}
\ee
with the term arising from $\psi_1(\eta_m,\eta_{m-1})$ disappearing in the limit.
Further, for terms with $m=m_j$ in the subsequence we have
$$
\|V_{m_j}\|=\delta_{m_j}\|Y_{m_j}\|\le\delta_{m_j}\|Y_{m_j}-V_{m_j}\|+\delta_{m_j}\|V_{m_j}\|
$$
and hence
\be
\|V_{m_j}\|\le\frac{\delta_{m_j}}{1-\delta_{m_j}}\|Y_{m_j}-V_{m_j}\|
\label{S52}
\ee
and so from \eqref{S51} and \eqref{S52} it follows that
\be
\lim_{j\to\infty}\frac{V_{m_j}}{\eta_{m_j}}=0,\qquad
\lim_{j\to\infty}\frac{Y_{m_j}}{\eta_{m_j}}=v.
\label{S53}
\ee
Now let $y_m=\langle v,Y_m\rangle$. Then taking 
the inner product with $v$ in the second limit in \eqref{S53}
gives ${\ds{\lim_{j\to\infty}}}\eta_{m_j}^{-1}y_{m_j}=1$, hence from \eqref{S54} that
\be
By_{m_j}>0\hbox{ for all large }j.
\label{S55}
\ee
Noting that $X_m\le p$ with $X_m\ne p$, hence $Y_m\le 0$ with $Y_m\ne 0$, for every $m$, and because $v_n>0$ for every $n$,
we have that $y_m<0$. Thus with \eqref{S55} we conclude that $B<0$.

We now repeat the above argument but instead taking $Y_m=SX_{-m}-p$. The only difference occurs at the end,
when we note that $Y_m\ge 0$ with $Y_m\ne 0$, and so $y_m>0$ for every $m$, which leads
to the conclusion that $B>0$. This is a contradiction, and
with this the theorem is proved.
\qed

\section{Genericity: The Proof of Theorem 1.2}

The purpose of this section is to prove Theorem 1.2.  Recall that we denote the set of normal families of bistable nonlinearities by $\mathcal{N}$.  Given $f_0 \in \mathcal{N}$, we consider $f\in\mathcal{N}$ of the form $f(u,a)=\gamma(u)f_0(u,a)$
where $\gamma\in C^{2}_+$, and we wish to show that Condition B holds for such $f$, for a residual set of $\gamma$
in $C^{2}_+$.
Below, Propositions 6.3 and 6.4 will establish that the set $\mathcal A$ of $\gamma$ for which Condition A holds
is a residual set. Then
Propositions 6.5 and 6.6 will establish that Condition B holds on a residual set of $\mathcal{A}$,
thereby completing the proof of Theorem 1.2.

The basic tool used here is the Abraham Transversality Theorem. If $Y\subseteq X$ is a closed subspace of a Banach
space $X$, we say $Y$ is complemented in $X$ if there exists a closed subspace $Z\subseteq X$ such that
$X=Y\oplus Z$. Certainly any subpace of either finite dimension or finite codimension is complemented, but
subspaces which are not complemented do exist. If $F:U\to W$ is a smooth map between
two Banach manifolds $U$ and $W$, and $M\subseteq W$ is a smooth submanifold of $W$ (with these manifolds
possibly infinite dimensional), then we say that $F$ is transverse to $M$ on a set $K\subseteq U$ if
whenever $F(x)\in M$ for some $x\in K$, then
$\ran(DF(x))+T_{F(x)}M=T_{F(x)}W$ (this sum of subspaces need not be a direct sum), and the space
$\Sigma_x=\{\overline x\in T_xU\;|\;DF(x)\overline x\in T_{F(x)}M\}$
is complemented in $T_xU$. 
(For purposes of this definition smooth means $C^1$.)
If $F$ is transverse to $M$ as above, then for some neighborhood $O\subseteq U$ with $K\subseteq O$,
the set $N=F^{-1}(M)\cap O$ is a submanifold of $U$, with $T_xN=\Sigma_x$ for every $x\in N$.

We also recall the Smale Density Theorem, which states that if $F:U\to W$ is a $C^r$ (for some $r\ge 1$) map between Banach
manifolds $U$ and $W$, with $U$ Lindel\"of, and if further $DF(x)$ is a Fredholm operator of index $j$
for every $x\in U$, then the set of regular values of $F$ is a residual subset of $W$ provided that $r>j$.
\vv

\noindent {\bf Theorem 6.1 (Abraham Transversality Theorem).} \bxx
Let $F:U\times V\to W$ be a $C^r$ map, where $U$, $V$, and $W$ are $C^r$ Banach manifolds. Assume also that
$U$ and $V$ are Lindel\"of spaces (for example, affine subspaces or open subsets of a separable Banach space). Suppose
that $M\subseteq W$ is a $C^r$ submanifold of $W$ and that $F$ is transverse to $M$ on $U\times V$.
Further suppose that for each $(x,\lambda)\in U\times V$ for which $F(x,\lambda)\in M$, the map
$$
\pi D_1F(x,\lambda)\in\LL(T_xU,\;T_{F(x,\lambda)}W/T_{F(x,\lambda)}M)
$$
is Fredholm of index $j$, and that $r>\max\{j,0\}$, where
$$
\pi\in\LL(T_{F(x,\lambda)}W,\;T_{F(x,\lambda)}W/T_{F(x,\lambda)}M)
$$
is the canonical projection onto the quotient space, and where $D_1F(x,\lambda)$ denotes the derivative
with respect to the first argument $x$. Let $F_\lambda:U\to W$ denote the map $F_\lambda(x)=F(x,\lambda)$
for each $\lambda\in V$. Then
$$
\{\lambda\in V\;|\;F_\lambda\hbox{ is transverse to }M\hbox{ on }U\}
$$
is a residual subset of $V$.
\exx
\vv

\noindent {\bf Sketch of Proof.}
The set $\Gamma=F^{-1}(M)$ is a submanifold of $U\times V$ with
$T_{(x,\lambda)}\Gamma=\{(\overline x,\overline \lambda)\in T_xU\times T_\lambda V\;|\;
\pi DF(x,\lambda)(\overline x,\overline\lambda)=0\}$.
Consider the map $\Pi:\Gamma\to V$ given by $\Pi(x,\lambda)=\lambda$. Then
one shows that $\lambda\in V$ is a regular value of $\Pi$ if and only if
$F_\lambda$ is transverse to $M$ on $U$. Further, the derivative
$D\Pi(x,\lambda)\in\LL(T_{(x,\lambda)}\Gamma,T_\lambda V)$ is Fredholm with
index $j$. By Smale's density theorem, the set of regular values $\lambda$
of $\Pi$ is a residual subset of $V$.
\qed
\vv

In what follows we shall let
$$
\begin{array}{lcl}
Q_0 &\bb = &\bb
\ds{\{q\in\ell^\infty(\Z)\;|\;\lim_{n\to\pm\infty}q_n=0\},}\\
\\
Q_1 &\bb = &\bb
\ds{\{q\in\ell^\infty(\Z)\;|\;\lim_{n\to\pm\infty}q_n=\pm 1\},}\\
\\
Q_2 &\bb = &\bb
\{(p,q)\in Q_1\times Q_1\;|\;p\ne S^kq\hbox{ for every }k\in\Z\}.
\end{array}
$$
Note that $Q_1$ is a Banach manifold which is Lindel\"of,
and that $T_qQ_1=Q_0$. The following result implies that
in a trivial fashion, $Q_2$ is also a Banach manifold which is Lindel\"of, and
that $T_{(p,q)}Q_2=Q_0\times Q_0$.
\vv

\noindent {\bf Lemma 6.2.} \bxx
The set $Q_2$ is an open subset of $Q_1\times Q_1$.
\exx
\vv

\noindent {\bf Proof.}
Take a sequence of points $(p_j,q_j)\in (Q_1\times Q_1)\setminus Q_2$
in the complement of $Q_2$, where $p_j=\{p_{j,n}\}_{n\in\Z}$
and $q_j=\{q_{j,n}\}_{n\in\Z}$. Assume these converge $(p_j,q_j)\to(p_0,q_0)$ to
some $(p_0,q_0)\in Q_1\times Q_1$. We must show that $(p_0,q_0)\not\in Q_2$.

We have that $p_j=S^{k_j}q_j$ for some sequence of integers $k_j$, from
the definition of  $Q_2$. By passing to a subsequence, we may assume that
either $k_j=k\in\Z$ is independent of $j$, or else $k_j\to\infty$ as $j\to\infty$,
or else $k_j\to-\infty$ as $j\to\infty$.

If $k_j=k$, then clearly $p_0=S^kq_0$,
and so $(p_0,q_0)\not\in Q_2$, as desired.
Thus suppose that $k_j\to\infty$
as $j\to\infty$. Then for any fixed $n$ we have
$$
\ba{lcl}
|p_{j,n}-q_{0,n+k_j}| &\bb \le &\bb
|p_{j,n}-q_{j,n+k_j}|+|q_{j,n+k_j}-q_{0,n+k_j}|\\
\\
&\bb \le &\bb
\|p_j-S^{k_j}q_j\|+\|S^{k_j}(q_j-q_0)\|\to 0
\ea
$$
as $j\to\infty$. But also $p_{j,n}-q_{0,n+k_j}\to p_{0,n}-1$ as $j\to\infty$ since $q_0\in Q_1$,
and so $p_{0,n}=1$. However, this contradicts the fact that $p_0\in Q_1$.

We omit the case in which $k_j\to-\infty$, as it is similar to the previous case.
With this, the proof is complete.
\qed
\vv

In what follows we shall let $D_a$ denote the derivative with respect
to the parameter $a$, with prime $'$ denoting the derivative with
respect to $u$ as usual, for the two arguments $(u,a)$ of $f$.
\vv

\noindent {\bf Proposition 6.3.} \bxx
Fix any $f \in \mathcal{N}$,
and define $G_f : Q_2 \times (-1,1) \to Q_1\times Q_1$ by 
$$
G_f(p,q,a) = \left(\ba{c}(S+S^{-1}-2I) p - f(p,a)\\ (S+S^{-1}-2I) q - f(q,a)\strt\ea\right).
$$
Suppose that $G_f$ is transverse to $\{(0,0)\}$ on $Q_2\times(-1,1)$.  Then Condition A holds.
\exx
\vv

\noindent {\bf Proof.}
We prove the contrapositive.  Suppose that Condition A fails.  Then there exists $(p,q) \in Q_2$
such that $G_f(p,q,a_+(0)) = (0,0)$.  It suffices to prove that $G_f$ is not transverse to $\{(0,0)\}$
at the point $(p,q,a_+(0))$, namely, that the derivative $DG_f(p,q,a_+(0))$ is not surjective.
We have that
\be
DG_f(p,q,a_+(0))(\overline{p},\overline{q},\overline{a}) =
\left(\ba{c}
\ds{L_p \overline{p} - D_af(p,a_+(0))\overline{a}}\\
\ds{L_q \overline{q} - D_af(q,a_+(0))\overline{a}}\strt\ea\right)
\label{dg}
\ee
where the operators $L_p=L_p(a_+(0))$ and $L_q=L_q(a_+(0))$ are as in \eqref{Ldef},
although here considered as an element of $\LL(Q_0)$ rather than $\LL(\ell^\infty(\Z))$.
By Lemma 3.1 both $L_p$ and $L_q$ are Fredholm with index zero and kernels of dimension either zero or one.
It follows from an implicit function theorem argument similar to that used in the proof of Proposition 1.4
that the dimension of the kernel of each of these operators is in fact one, and thus
$$
\codim\ran (L_p)=\codim\ran (L_q)=1
$$
for the codimensions of the ranges in the space $Q_0$.
From this and the formula \eqref{dg}, it follows immediately that the range of
$DG_f(p,q,a_+(0))$ has codimension either one or two. Thus $DG_f(p,q,a_+(0))$ is
not surjective, as desired.
\qed
\vv

Recall now the set $C^{2}_+$ given by \eqref{C3}, which is endowed with the usual $C^{2}$
topology.
\vv

\noindent {\bf Proposition 6.4.} \bxx
Fix any $f_0\in\mathcal{N}$, and define
$\tilde G_{f_0} : Q_2 \times (-1,1) \times C^{2}_+ \to Q_1\times Q_1$ by
$$
\tilde G_{f_0}(p,q,a,\gamma) = \left(\ba{c}
(S+S^{-1}-2I) p - \gamma(p)f_0(p,a)\\
(S+S^{-1}-2I) q - \gamma(q)f_0(q,a)\strt\ea\right).
$$
Then
$\tilde G_{f_0}$ is transverse to $\{(0,0)\}$ on $Q_2\times(-1,1)\times C^{2}_+$.
\exx
\vv

\noindent {\bf Proof.}
Let us denote the dependence of the operator $L_q=L_q(a)$ in \eqref{Ldef} on the nonlinearity $f$ by
\be
L_{f,q}=L_{f,q}(a)=S+S^{-1}-2I-f'(q,a),
\label{lfq}
\ee
for $f\in\mathcal{N}$. Now take any $(p,q,a,\gamma)\in Q_2\times(-1,1)\times C^{2}_+$
for which $\tilde G_{f_0}(p,q,a,\gamma)=0$. Denoting $f(u,a)=\gamma(u)f_0(u,a)$, we have that
\be
D\tilde G_{f_0}(p,q,a,\gamma)(\overline{p},\overline{q},\overline{a},\overline{\gamma}) = 
\left(\ba{c}
L_{f,p}\overline{p} - D_af(p,a)\overline{a} - \overline{\gamma}(p) f_0(p,a)\\
L_{f,q}\overline{q} - D_af(q,a)\overline{a} - \overline{\gamma}(q) f_0(q,a)\strt\ea\right).
\label{dg2}
\ee
By Lemma 3.1 both $L_{f,p}$ and $L_{f,q}$ are Fredholm with index zero and kernels of dimension either zero or one,
considered here as elements of $\LL(Q_0)$. We must prove that the operator $D\tilde G_{f_0}(p,q,a,\gamma)$ is surjective.
Denote by $v_p$ any nonzero element of the kernel $\ker (L_{f,p})$, if such exists, and similarly with $v_q$ as
a nontrivial element of $\ker (L_{f,q})$.

Three cases now arise. In the first case we suppose that both $L_{f,p}$ and $L_{f,q}$ are isomorphisms,
and so surjectivity is immediate.

In the second case we assume exactly one of these operators is an isomorphism, say $L_{f,q}$ for definiteness,
and so the operator $L_{f,p}$ has a nontrivial kernel element $v_p$. Moreover, the range of $L_{f,p}$ is characterized as
$$
\ran (L_{f,p})=\{w\in Q_0\;|\;\langle v_p,w\rangle=0\}
$$
by Lemma 3.1.
If it is not the case that $D\tilde G_{f_0}(p,q,a,\gamma)$ is surjective
then the first line in \eqref{dg2} is annihilated by $v_p$, for every
choice of $(\overline p,\overline a, \overline \gamma)$. In particular,
taking $\overline p=0$ and $\overline a=0$, and $\overline \gamma$ to be a function for which
\be
\overline\gamma(p_n)=\left\{\ba{cl}
0, & n\ne k,\\
1, & n=k,\strt
\ea\right.
\label{gamk}
\ee
for some integer $k$, we have that
\be
0=\langle v_p,\overline\gamma(p)f_0(p,a)\rangle=v_{p,k}f_0(p_k,a).
\label{pk}
\ee
But $f_0(p_k,a)\ne 0$ for all large $k$, and so \eqref{pk} implies that $v_{p,k}=0$
for all large $k$. This contradicts Lemma 3.1, which states that $v_{p,k}\ne 0$
for all large $k$.

It remains only to consider the case where neither $L_{f,p}$ nor $L_{f,q}$ are isomorphisms, and so the vectors
$v_p$ and $v_q$, respectively, annihilate their ranges of these operators. Then if the derivative
$D\tilde G_{f_0}(p,q,a,\gamma)$ is not surjective, its range is annihilated by some nontrivial nontrivial
combination of $v_p$ and $v_q$. To be precise, there exist constants $\tau_p$ and $\tau_q$, not both zero, such that
$$
\ba{lcl}
0 &\bb = &\bb \tau_p\langle v_p,\;L_{f,p}\overline p-D_af(p,a)\overline{a}-\overline{\gamma}(p) f_0(p,a)\rangle\\
\\
&\bb &\bb +\tau_q\langle v_q,\;L_{f,q}\overline q-D_af(q,a)\overline{a}-\overline{\gamma}(q) f_0(q,a)\rangle
\ea
$$
for every $(\overline p,\overline q,\overline a, \overline \gamma)$.
Taking $\overline p=\overline q=0$ and $\overline a=0$ gives
$$
0=\tau_p\sum_{n=-\infty}^\infty v_{p,n}\overline\gamma(p_n)f_0(p_n,a)
+\tau_q\sum_{n=-\infty}^\infty v_{q,n}\overline\gamma(q_n)f_0(q_n,a)
$$
for every $\overline\gamma$.
If one of the coefficients $\tau_p$ or $\tau_q$ is zero, say $\tau_q=0$,
then choosing a function $\overline\gamma$ as in \eqref{gamk}
and arguing as before yields a contradiction. Thus assume that both $\tau_p\ne 0$
and $\tau_q\ne 0$. If there exist arbitrarily large $k$ for which the point
$p_k$ is distinct from all the points $q_j$, that is, $p_k\ne q_j$ for every $j\in\Z$,
then we may choose $\overline\gamma$ satisfying \eqref{gamk} and additionally satisfying
$\overline\gamma(q_j)=0$ for every $j\in\Z$. This now gives equation \eqref{pk},
and so $v_{p,k}=0$ for such $k$. However, this contradicts Lemma 3.1.

Thus we have that for every sufficiently large integer $k$, there exists
an integer $j$, such that $p_k=q_j$. Necessarily $j=j_k$ is uniquely determined
from $k$, and $j_k\to\infty$ as $k\to\infty$. Reversing the roles of $p$ and $q$
imply that also, for every sufficiently large integer $j$, there exists
an integer $k=k_j$ such that $q_j=p_{k_j}$. It follows immediately that there exists
some integer $m$ such that $p_k=q_{k+m}$ for all sufficiently large $k$.
As both $p$ and $q$ satisfy the difference equation \eqref{scalardiff}, we conclude
that $p_k=q_{k+m}$ for every $k\in\Z$, that is, $p=S^mq$. But this contradicts
the fact that $(p,q)\in Q_2$, and completes the proof.
\qed
\vv

In what follows, let us denote
$$
\ba{lcl}
\mathrm{Fred}_{1,1} &\bb = &\bb
\{\Omega\in\LL(Q_0)\;|\;\Omega\hbox{ is a Fredholm operator}\\
\\
&\bb &\bb \hbox{of index zero, with }\dim\ker(\Omega)=1\}.
\ea
$$
It is known that $\mathrm{Fred}_{1,1}$ is a $C^\infty$ submanifold of $\LL(Q_0)$ of codimension one.
Moreover, if $\Omega\in\mathrm{Fred}_{1,1}$ then the tangent space of $\mathrm{Fred}_{1,1}$ at this point
is given by
$$
T_\Omega\mathrm{Fred}_{1,1}=\{\overline\Omega\;|\;\overline\Omega v\in\ran(\Omega)
\hbox{ whenever }v\in\ker(\Omega)\}.
$$
\vv

\noindent {\bf Proposition 6.5.} \bxx
\label{fredlemma}
Fix any $f \in \mathcal{N}$, and suppose that Condition A holds for $f$.
Define $F_f : Q_1 \times (-1,1) \to Q_1\times \LL(Q_0)$ by 
$$
F_f(p,a) = \left(\ba{c}(S+S^{-1}-2I) p - f(p,a)\\ S+S^{-1}-2I - f'(p,a)\strt\ea\right).
$$
Also define the manifold
$$
M=\{0\}\times\mathrm{Fred}_{1,1}\subseteq Q_1\times\LL(Q_0).
$$
Then $f$ satisfies Condition B if and only if
$F_f$ is transverse to $M$ at the point $(p,a_+(0))$, where $p$ is as in the
statement of Condition A\@.
\exx
\vv

%

\noindent {\bf Proof.}
In what follows $p$ is as in Condition A, with $L_p=L_p(a_+(0))$, so $L_p\in\mathrm{Fred}_{1,1}$.
Also, let $v\in Q_0$ be the kernel element as in the statement of Proposition 1.4. Then
one sees immediately that
$$
T_{L_p} \mathrm{Fred}_{1,1}=\{\overline\Omega\in\LL(Q_0)\;|\;\langle v,\overline\Omega v\rangle=0\},\qquad
L_p=S+S^{-1}-2I - f'(p,a_+(0)).
$$
We now compute
$$
DF_f(p,a_+(0))(\overline{p},\overline{a}) =
\left(\ba{c}
L_p \overline{p} - D_af(p,a_+(0))\overline{a}\\
-f''(p,a_+(0))\overline{p} - D_af'(p,a_+(0))\overline{a}\strt
\ea
\right).
$$
Thus $F_f$ is transverse to $M$ at $(p,a_+(0))$ if and only if,
for every choice of $w \in Q_0$ and $\Omega \in \LL(Q_0)$, there exists $(\overline{p},\overline{a})\in Q_0\times\R$,
and $\overline{\Omega} \in T_{L_p} \mathrm{Fred}_{1,1}$, such that 
\begin{equation}
\ba{lcl}
w &\bb = &\bb L_p \overline{p} - D_af(p,a_+(0))\overline a, \\
\\
\Omega &\bb = &\bb -f''(p,a_+(0))\overline{p} - D_af'(p,a_+(0))\overline a + \overline{\Omega}.
\ea
\label{eq:transdefB}
\end{equation}
Note that the existence of $\overline\Omega$ satisfying the second line in \eqref{eq:transdefB}
is equivalent to the equation
\be
\langle v,\Omega v\rangle+\langle v,f''(p,a_+(0))(\overline p,v)\rangle
+\langle v,D_af'(p,a_+(0))(\overline a,v)\rangle=0.
\label{fred}
\ee
Further note that we have
$$
\langle v,f''(p,a_+(0))(v,v)\rangle=2B,
$$
where $B$ is the quantity \eqref{bsum} in Condition B\@.

To prove the lemma, suppose first that Condition B holds, and let $w$ and $\Omega$
be given as above. Then the first equation in \eqref{eq:transdefB} can be solved
for some $\overline p$ and $\overline a$ if and only if
$\langle v,w\rangle=-\langle v,D_af(p,a_+(0))\rangle\overline a$, equivalently,
\be
\overline a=-\frac{\langle v,w\rangle}{\langle v,D_af(p,a_+(0))\rangle}.
\label{ola}
\ee
Here we have used the fact that $f$ is a normal family, in particular that
$D_af(u,a) > 0$ in $(-1,1)\times(-1,1)$, and also
Proposition 1.4, in particular that $v_n > 0$,
to conclude that the denominator is nonzero.
With this choice of $\overline a$ there exists $\overline p$ satisfying the
first equation in \eqref{eq:transdefB}, and moreover, $\overline p$ is uniquely
determined up to an additive multiple of $v$. Thus fixing a particular choice $\tilde p$ of
$\overline p$, we see that the general form of $\overline p$ is
$\overline p=\tilde p+\lambda v$ for arbitrary $\lambda\in\R$. Thus
equation \eqref{fred} takes the form
\be
\langle v,\Omega v\rangle+\langle v,f''(p,a_+(0))(\tilde p,v)\rangle+2\lambda B
+\langle v,D_af'(p,a_+(0))(\overline a,v)\rangle=0
\label{fred2}
\ee
with the above choices of $\overline a$ and $\overline p$. As $B\ne 0$ is assumed,
there is a unique choice of $\lambda$ solving this equation, as desired.

Now assume that Condition B does not hold, and so $B=0$. Taking $w=0$ above
forces $\overline a=0$, by \eqref{ola}, and so necessarily $\overline p=\lambda v$
to satisfy the first equation in \eqref{eq:transdefB}. Equation \eqref{fred2}
takes the form $\langle v,\Omega v\rangle=0$. However, this equation is not
satisfied for every $\Omega$, which implies the transversality condition fails,
as claimed.
\qed
\vv


\noindent {\bf Proposition 6.6.} \bxx
Fix any $f_0\in\mathcal{N}$, and define
$\tilde F_{f_0} : Q_1 \times (-1,1) \times C^{2}_+ \to Q_1\times \LL(Q_0)$ by
$$
\tilde F_{f_0}(p,a,\gamma) = \left(\ba{c}
(S+S^{-1}-2I) p - \gamma(p)f_0(p,a)\\
S+S^{-1}-2I - (\gamma(p)f_0(p,a))'
\strt
\ea\right).
$$
Then
$\tilde F_{f_0}$ is transverse to $M$ on $Q_1\times(-1,1)\times C^{2}_+$,
where $M$ is as in the statement of Proposition 6.5.
\exx
\vv


\noindent {\bf Proof.}
As before, let $L_{f,q}=L_{f,q}(a)$ denote the operator \eqref{lfq}. Then,
denoting $f(u,a)=\gamma(u)f_0(u,a)$ in what follows, we compute the derivative
$$
D\tilde F_{f_0}(p,a,\gamma)(\overline{p},\overline{a},\overline{\gamma}) = 
\left(\ba{c}
L_{f,p}\overline{p} - \gamma(p)D_af_0(p,a)\overline a - \overline{\gamma}(p)f_0(p,a)\\
-f''(p,a)\overline p -D_af'(p,a)\overline a -(\overline\gamma(p)f_0(p,a))'
\strt
\ea\right).
$$
Thus $\tilde F_{f_0}$ is transverse to $M$ at a point $(p,a,\gamma)$ for which $\tilde F_{f_0}(p,a,\gamma)\in M$
if and only if,
for every choice of $w \in Q_0$ and $\Omega \in \LL(Q_0)$, there exist
$(\overline{p},\overline{a},\overline \gamma)\in Q_0\times\R\times C^{2}[-1,1]$,
such that 
\be
\ba{lcl}
\makebox[0pt][l]{$w=L_{f,p}\overline{p} - \gamma(p)D_af_0(p,a)\overline a - \overline{\gamma}(p)f_0(p,a),$} \\
\\
\langle v,\Omega v\rangle+\langle v,f''(p,a)(\overline p,v)\rangle
&\bb + &\bb
\langle v,D_af'(p,a)(\overline a,v)\rangle\\
\\
&\bb + &\bb \langle v,(\overline\gamma(p)f_0(p,a))'v\rangle=0,
\ea
\label{womega}
\ee
much as in \eqref{eq:transdefB} and \eqref{fred}. Here $v$ is a nontrivial element of the kernel
of $L_{f,p}$, and is
is uniquely determined up to scalar multiple (which we fix)
as $L_{f,p}=S+S^{-1}-2I-f'(p,a)\in\mathrm{Fred}_{1,1}$. We note that there is no assurance
that $v_n>0$ for the coordinates.

Given any such $w$ and $\Omega$, first choose $\overline a=0$. Next,
fix $k$ large enough that both $v_k\ne 0$ (which can be done by Lemma 3.1)
and $f(p_k,a)\ne 0$,
and choose $\overline\gamma\in C^{2}[-1,1]$ to satisfy
$$
\overline\gamma(p_n)=\left\{\ba{cl}
0, & n\ne k,\\
\gamma^0, & n=k,\strt
\ea\right.
\qquad
\overline\gamma'(p_n)=\left\{\ba{cl}
0, & n\ne k,\\
\gamma^1, & n=k,\strt
\ea\right.
$$
with the quantities $\gamma^0$ and $\gamma^1$ to be chosen shortly.
There exists $\overline p$ satisfying the first equation in \eqref{womega}
if and only if $\langle v,w\rangle=-\langle v,\overline\gamma(p)f_0(p,a)\rangle$, or equivalently,
$$
\gamma^0=-\frac{\langle v,w\rangle}{v_kf_0(p_k,a)}.
$$
Taking this value for $\gamma^0$,
we now fix the choice of $\overline p$ as well. In particular, we note that although
we have not yet chosen $\gamma^1$, the choice of this quantity will not affect $\overline p$
as the first equation in \eqref{womega} does not involve $\overline\gamma'(p)$.

Finally, we make the unique choice of $\gamma^1$ so that the second
equation in \eqref{womega} holds, and one sees that the required value is
$$
\gamma^1=-\frac{\langle v,\Omega v\rangle+\langle v,f''(p,a)(\overline p,v)\rangle
+v_k^2\gamma^0f'_0(p_k,a)}{v_k^2f_0(p_k,a)}.
$$
With this, the proof is complete.
\qed
\vv

\noindent {\bf Proof of Theorem 1.2.}
With $f_0\in\mathcal{N}$, we apply the Abraham Transversality Theorem to the maps $\tilde G_{f_0}$
and $\tilde F_{f_0}$ in Propositions 6.4 and 6.6, respectively. What must be checked is that the appropriate
operators are Fredholm with the appropriate index. For the map $\tilde G_{f_0}$ we consider the derivative
$D_{1,2,3}\tilde G_{f_0}(p,q,a,\gamma)$ taken with respect to the first three arguments
$p$, $q$, and $a$, but not $\gamma$ (as $\gamma$ plays the role of $\lambda$ in the above statement
of the Transversality Theorem).
As noted in the proof of Proposition 6.4,
the operators $L_{f,p}$ and $L_{f,q}$ are Fredholm with index zero, and from this it follows
easily, using the formula \eqref{dg2}, that $D_{1,2,3}\tilde G_{f_0}(p,q,a,\gamma)$
is Fredholm of index one. Thus the Transversality Theorem applies as the map $\tilde G_{f_0}$
is $C^2$.

For the map $\tilde F_{f_0}$ we must follow the derivative $D_{1,2}\tilde F_{f_0}(p,a,\gamma)$
taken with respect to $p$ and $a$, with the projection which annihilates the tangent space
$T_{L_p}\mathrm{Fred}_{1,1}$. This gives the operator
\be
\pi D_{1,2}\tilde F_{f_0}(p,a,\gamma)(\overline{p},\overline{a}) = 
\left(\ba{c}
L_{f,p}\overline{p} - \gamma(p)D_af_0(p,a)\overline a\\
-\langle v,f''(p,a)(\overline p,v) +D_af'(p,a)(\overline a,v)\rangle
\strt
\ea\right),
\label{final}
\ee
where $v$ is the kernel element of $L_{f,p}$. Again $L_{f,p}$ has index zero, and
as the second coordinate in the range of \eqref{final} is scalar, it follows that the operator
\eqref{final} has Fredholm index zero. Again the Transversality Theorem applies as the
map $F_{f_0}$ is $C^1$.
It follows directly that for a residual set of $\gamma\in C^{2}_+$ the maps
$G_{\gamma f_0}$ and $F_{\gamma f_0}$ are transverse to $\{(0,0)\}$ and to $M$,
respectively, on their domains, and thus that Condition B holds for $\gamma f_0$.
With this Theorem 1.2 is proved.
\qed

\bibliography{CPgenBib}{}

\end{document}